\documentclass[12pt]{article} 
\usepackage[notref,notcite]{showkeys}
\usepackage{amsfonts,amsmath,layout}
\usepackage{graphicx}   
\usepackage{url}
\usepackage{hyperref}

\newcommand{\N}{\mathbb N}

\newcommand{\R}{\mathbb R}

\def\QED{$\; \Box$\medskip}

\def\R{\mathbb R}
\def\N{\mathbb N}

\def\E{\mathbb E}
\def\P{\mathbb P}
\def\O{\mathcal O}

\def\qqf{\qquad\forall}

\def\qmb{\qquad\mbox}
\def\ds{\displaystyle}
\newtheorem{Theorem}{Theorem}[section]

\newtheorem{Proposition}[Theorem]{Proposition}
\newtheorem{Lemma}[Theorem]{Lemma}
\newtheorem{Corollary}[Theorem]{Corollary}
\newtheorem{Remark}[Theorem]{Remark}
\newtheorem{Example}[Theorem]{Example}

\begin{document}

\title{H\"{o}lder estimates in space-time for viscosity solutions \\of Hamilton-Jacobi equations
\thanks{This work was partially supported by the French ANR (Agence Nationale de la Recherche) through MICA project (ANR-06-06-BLAN-0082) and KAMFAIBLE project (BLAN07-3-187245), as well as by the Italian PRIN 2005 and PRIN 2007 Programs ``Metodi di viscosit\`a, metrici e di teoria del controllo in equazioni alle derivate parziali nonlineari''..}}
\author{Piermarco Cannarsa\thanks{Dipartimento di Matematica, Universit\`a di Roma ``Tor Vergata'',
Via della Ricerca Scientifica 1, 00133 Roma (Italy),
e-mail: cannarsa@mat.uniroma2.it}\quad\&\quad Pierre Cardaliaguet \thanks{Universit\'e de Bretagne
Occidentale, UMR 6205, 6 Av. Le Gorgeu,
BP 809, 29285 Brest (France); e-mail:
Pierre.Cardaliaguet@univ-brest.fr}}

\maketitle

\begin{abstract} 
It is well-known that solutions to the basic problem in the calculus of variations
may fail to be Lipschitz continuous when the Lagrangian depends on $t$. Similarly, for viscosity solutions to time-dependent Hamilton-Jacobi  equations one cannot expect Lipschitz bounds to hold uniformly with respect to the regularity of coefficients. This phenomenon raises the question whether such solutions satisfy uniform estimates in some weaker norm. 

We will show that this is the case for a suitable H\"older norm, obtaining uniform  estimates in $(x,t)$ for solutions to first and second order Hamilton-Jacobi equations. Our results apply to degenerate parabolic equations and require superlinear growth at infinity, in the gradient variables, of the Hamiltonian.  Proofs are based on comparison arguments and representation formulas for viscosity solutions, as well as weak reverse H\"older inequalities.
\end{abstract} 



 \section{Introduction}\label{se:intro}
The object of this paper is the regularity of solutions to the Hamilton-Jacobi equation
\begin{equation}\label{intro:HJ2}
u_t(x,t)-{\rm Tr}\left(a(x,t)D^2u(x,t)\right)+H(x,t,Du(x,t))=0\qquad {\rm in }\; \R^N\times (0,T)
\end{equation}
where 
$H$ and  $a$ will be assumed to satisfy the following hypotheses:
\begin{itemize}
\item  there are  real numbers $q>2$, $\delta > 1$ 
and $\eta_\pm\geq 0$ such that
\begin{equation}\label{intro:GrowthCond2}
\frac{1}{\delta}|z|^q-\eta_{_-} \leq H(x,t,z) \leq \delta|z|^q+\eta_{_+}\qquad \forall (x,t,z)\in \R^N\times(0,T)\times \R^N;
\end{equation}
\item $a=\sigma\sigma^*$ for some locally Lipschitz continuous map $\sigma:(x,t)\mapsto\sigma(x,t)$, with values in the $N\times D$ real matrices ($D\geq 1$), such that $\|\sigma(x,t)\|\le \delta$
for all $(x,t)\in\R^N\times(0,T)$.
\end{itemize}
We note that no initial condition is needed for our analysis, nor convexity of  $H$ in $Du$.

For a given viscosity solution $u$ of \eqref{intro:HJ2}, the kind of regularity properties we are interested in are uniform continuity estimates in $(x,t)$   that do not depend on the smoothness of coefficients but just on the constants that appear in  \eqref{intro:GrowthCond2}, and on the sup-norm of $u$. Another important feature of our approach is that the above equation will not be assumed to be uniformly parabolic. Indeed, it will be allowed to degenerate to the point of reducing to the first order equation 
\begin{equation}\label{intro:HJ0}
u_t+H(x,t,Du)=0\qquad {\rm in }\quad \R^N\times (0,T)\,,
\end{equation}
in which case we will just require $q>1$ in \eqref{intro:GrowthCond2}.

The typical form of our results ensures that any bounded continuous viscosity solution $u$ of \eqref{intro:HJ2}  satisfies, for positive time, the uniform H\"older estimate 
\begin{equation}\label{intro:main2}
|u(x,t)-u(y,s)|\leq C\left[|x-y|^{\theta-p\over \theta-1}+|t-s|^{\theta-p\over\theta}\right]
\end{equation}
where $p$ is the conjugate exponent of $q$, and $\theta>p$ depends only on the aforementioned constants. The above result may take a specific form according to the problem we will consider. For instance, for second order equations we suppose that the Hamiltonian $H$ is super-quadratic ($q>2$), whereas for problem \eqref{intro:HJ0} we just need super-linear growth ($q>1$). Moreover, for both first and second order problems we can also give a local version of our result, that is, an estimate that applies to solutions in an open set $\O\subset \R_x^N\times \R_t$.

In order to better understand the problem under investigation it is convenient to start the analysis with first order equations. In this case, when $H(x,t,z)$ is convex in $z$, the viscosity solutions of \eqref{intro:HJ0} can be represented as value functions of problems in the calculus of variations.
Consequently, the regularity of $u$ is connected with that of minimizers. As is well-known, minimizers are Lipschitz continuous in the autonomous case (see \cite{clvi85}, \cite{AAB}, \cite{dafr03}, \cite{da07}), so that solutions turn out to be locally Lipschitz when $H=H(x,z)$. On the other hand, for nonautonomous problems, the Lipschitz regularity of minimizers is no longer true as is shown in \cite{AA}, and \eqref{intro:main2} is the optimal  H\"older estimate that can be expected, see  the example in section~\ref{se:example_Lip} and Remark~\ref{re:optimal_Holder} of this paper. A class of nonautonomous first order problems for which such an estimate can be obtained is studied in \cite{ca08}. Unlike the above references, however, our present results do not require $H$ to be convex in $z$. 

As for second order problems, H\"older regularity results for solutions  of {\em uniformly} parabolic equations have been the object of a huge literature for both linear and nonlinear problems. However, very few results can be found in connection with the present context, where we drop uniform parabolicity and allow for unbounded Hamiltonians. In the {\em stationary} case, Lipschitz bounds for solutions of uniformly elliptic equations with a super-quadratic Hamiltonian were obtained in \cite{LL}. More recently, H\"older estimates have been proved in \cite{CDLP} for viscosity subsolutions of fully nonlinear degenerate elliptic equations with super-quadratic growth in the gradient.

Our work is mainly motivated by homogenization theory, where such uniform estimates are necessary to study the limiting behavior of solutions and/or to prove the existence of correctors (see, e.g., \cite{ma04} and \cite{Sc}). For instance, estimate \eqref{intro:main2} could be applied to equations of the form
\begin{equation*}
u_t^\epsilon(x,t)-{\rm Tr}\left(a\Big(x,t, {x\over\epsilon}, {t\over\epsilon^2}\Big)D^2u^\epsilon(x,t)\right)+H\Big(x,t,{x\over\epsilon}, {t\over\epsilon^2},Du^\epsilon(x,t)\Big)=0
\end{equation*}
where $a(x,t,\cdot,\cdot)$ and $H(x,t,\cdot,\cdot,z)$ are periodic in $\R^N\times\R$.

A brief comment of the structure of the proof is now in order. Our reasoning involves three main steps: 
\begin{enumerate}
\item construction of suitable arcs along which super-solutions  exhibit a sort of monotone behavior; 
\item one-sided H\"older bound for sub-solutions;
\item application of a weak reverse H\"older inequality result.
\end{enumerate}
Let us be more specific on the above points in the simpler case of first order equations.
Our first step consists in showing that, if $u$ is a super-solution of \eqref{intro:HJ2},
then for any point $(\bar x,\bar t)\in\R^N\times(0,T]$ there is an arc $\xi\in W^{1,p}([0,\bar t];\R^N)$, satisfying $\xi(\bar t)=\bar x$, such that
\begin{equation}\label{intro:toto1}
u(\bar x,\bar t)\geq u(\xi(t),t) +C\int_t^{\bar t}|\xi'(s)|^pds -\eta_{_+} (\bar t-t)\qquad \forall t\in [0,\bar t]
\end{equation}
for some constant $C>0$. Second, using Hopf's formula, we obtain the following one-sided bound for  any sub-solution $u$ of \eqref{intro:HJ2}:
\begin{equation}\label{eq:one-sided}
u(\bar x,s)\;\leq\; u(y,t)+C(s-t)^{1-p}|y-\bar x|^p+\eta_{_-}(s-t) \qqf y\in \R^N\;, \forall s> t\,.
\end{equation}
So, choosing $y=\xi(t)$ in \eqref{eq:one-sided} and combining such an estimate with
\eqref{intro:toto1}, we derive 
\begin{equation*}
\frac{1}{\bar t-t} \int_t^{\bar t} |\xi'(s)|^pds \leq C_1\left(\frac{|\xi(t)-x|}{\bar t-t}\right)^p+C_0
\qquad \forall t\in [0,\bar t)
\end{equation*}
which yields, in turn, the weak reverse H\"older inequality
\begin{equation}\label{intro:toto2}
\frac{1}{\bar t-t} \int_t^{\bar t} |\xi'(s)|^pds \leq C_1\left(\frac{1}{\bar t-t}\int_t^{\bar t}|\xi'(s)|ds\right)^p+C_0
\qquad \forall t\in [0,\bar t)\,.
\end{equation}
Observe that \eqref{intro:toto2} is weaker than the classical reverse H\"older inequality 
used to improve the integrability of  functions  (see, e.g., \cite{gi83}). 
Nevertheless, we prove that,  
\begin{equation*}
%\label{co_toto2}
\int_t^{\bar t} |\xi'(s)|ds \leq C \ (\bar t-t)^{1-{1\over\theta}}
\qquad \forall t\in [0,\bar t]
\end{equation*}
for some exponent $\theta>p$ depending only on structural constants. Finally, we show that
the above inequalities imply  estimate \eqref{intro:main2} with exactly the same exponent $\theta$.

One of the interesting aspects of our approach is that, using essentially the same ideas we have just described, we manage to study the second order problem \eqref{intro:HJ2}. 
As it should be clear from the above discussion, such a transposition requires a certain familiarity with some techniques that are typical of stochastic analysis.
 For instance, the role of $\xi$ will be now played by the controlled diffusion process which satisfies
$dX_t= \zeta_tdt+\sigma(X_t, t)dW_t$, where $W$ is a standard $N$ dimensional Brownian motion and $\zeta$ is a $p$-summable adapted control. Moreover, the one-sided H\"older bound of step 2 will be recovered by the use of a suitable Brownian bridge. Furthermore, the stochastic version of our reverse H\"older inequality result will require
$$
\E\left[ \frac{1}{t-\bar t} \int_{\bar t}^t |\zeta_s|^pds\right] \leq C_1\,\E\left[\Big(\frac{1}{t-\bar t}\int_{\bar t}^t|\zeta_s|ds\Big)^p\right]+\frac{C_0}{(t-\bar t)^\frac{p}{2}}\qquad \forall t\in (\bar t,T]
$$
to yield the conclusion that
$$
\E \left[\Big(\int_{\bar t}^t|\zeta_s|ds\Big)^p\right] \leq C (t-\bar t)^{p-\frac{p}{\theta}}\left(\|\zeta\|^p_p
+B\right)
\qquad \forall t\in (\bar t,T]
$$
for some  $\theta\in(p,2)$.

The outline of this paper is the following. In section~\ref{se:preli} we fix notation and recall preliminaries from stochastic analysis, including the basic properties of Brownian bridges. Section~\ref{se:holder} is devoted to weak reverse H\"older inequalities. Then, we present our main results: we study the H\"older continuity of solutions to first order equations in section~\ref{se:first}, while second order problems are investigated in section~\ref{se:second} (for both problems we give a global and a local version of our results).  In between (section~\ref{se:exa}), we discuss counterexamples to higher regularity.

 \section{Notation and preliminaries}\label{se:preli}
We denote by 
 $x\cdot y$ the Euclidean scalar product  of two vectors $x,y\in\R^N$ and by
$|x|$ the Euclidean norm of $x$. For any $x_0\in\R^N$ and $r>0$, we denote by $B(x_0,r)$ the open ball of radius $r$, centered at $x_0\in\R^n$, and we set $B_r=B(0,r)$.

Let $D\ge 1$ be an integer. We denote by $\R^{N\times D}$ the space of all $N\times D$ real matrices equipped with the following norm
$$\|\sigma\|=\sqrt{\mbox{Tr}(\sigma \sigma^*) }\,,$$
where $\sigma^*$ denotes the transpose of $\sigma$ and 
$\mbox{Tr}(A)$ the trace of  $A\in\R^{N\times N}$.

We denote by $C(\R^N\times [0,T])$ the space of all continuous functions $u:\R^N\times [0,T]\to\R$.
 
For any  nonempty set $S\subset\R^N$ let $S^c=\R^N\setminus S$.  We denote $d_S$ the Euclidean {\em distance} function from $S$, that is, 
\begin{equation*}
d_S(x)=\inf_{y\in S}|x-y|\qquad\forall x\in\R^N\,.
\end{equation*}
It is well-known that $d_S$ is a Lipschitz function of constant 1. 

%{\bf Burkholder-Davis-Gundy inequality } (Voir par exemple Karatzas-Schreeve) states that if $(M_t)$ is a continuous martigale  in $\R^N$, then there are constants $0<c<C$ which depend only on $N$ 
%such that
%$$
%c\E\left(<M>_t^{p/2}\right)\leq \E\left((M_t^*)^{p}\right)\leq C \E\left(<M>_t^{p/2}\right) \qquad \forall p>1
%$$
%where $<M>_t$ is the quadratic variation of $M_t$ and 
%$$
%M^*_t= \sup_{s\in [0,t]}|M_s|
%$$
%In particular, let $(B_s)$ be a $P-$dimensional brownian motion and $f:\Omega\times [0,+\infty)\to \R^{D\times N}$ be adapted to the filtration
%generated by $(B_t)$. If $M_t=\int_0^t f_sdB_s$, then 
%$$
%<M>_t= \int_0^t {\rm Tr}(f_sf_s^*)ds\qquad \forall t\geq 0\;.
%$$
For $1\le p<\infty$ we denote by $L^p(a,b;\R^N)$ the space of all $p$-summable (with respect to the Lebesgue measure) Borel vector-valued functions $\xi:[a,b]\to\R^N$, and we use the shorter notation $L^p(a,b)$ if $N=1$. Similarly, we denote by $W^{1,p}([a,b];\R^N)$ the Sobolev space of all absolutely continuous arcs $\xi:[a,b]\to\R^N$ such that $\dot\xi\in L^p(a,b;\R^N)$ . 

Let now $(\Omega,{\cal F}, \P)$ be a stochastic basis, i.e., a measure space where $\P$ is a probability measure. We denote by  $L^p(\Omega\times[a,b];\R^N)$ the space of all
measurable  functions (with respect to the product measure) $\xi:\Omega\times[a,b]\to\R^N$, again suppressing the arrival set when $N=1$. 
In all the above cases, we denote by $\|\xi\|_p$ the standard $L^p$-norm of $\xi$.

Let $({\cal F}_t)$ be a filtration on $\Omega$. We denote by $L^p_{\rm ad}(\Omega\times[a,b];\R^N)$ the space of $p$-summable stochastic processes, adapted to $({\cal F}_t)$.

We will repeatedly use, in the sequel, the following classical estimate for solutions of the stochastic differential equation
\begin{equation}\label{eq:bridge}
dY_t= \zeta_tdt+\sigma(Y_t, t)dW_t\,,
\end{equation}
where $(W_t)$ is a $D$-dimensional Brownian motion adapted to
$({\cal F}_t)$. 
\begin{Lemma}\label{stoch_bound} Let  $\sigma:\R^N\times [0,T]\to\R^{N\times D}$ be a Lipschitz continuous map such that $\|\sigma\|\le \delta$, let $\zeta\in L^p_{\rm ad}(\Omega\times [0,T]; \R^N)$ $(p>1)$, and let $Y$ be a solution of \eqref{eq:bridge}.
Then, for every $r\in (0,p]$ there is a positive constant $C(r)$ such that
\begin{equation}\label{eq:estibridge}
\E\left[|Y_t-Y_s|^r\right]\le C(r)
\left\{\E\left[\Big|\int_s^t \zeta_\tau\,d\tau\Big|^r\right]+\delta^r\,|t-s|^{\frac r2}
\right\}\qqf s,t\in [0,T]\,.
\end{equation}
\end{Lemma}
{\it Proof:} For every $r\in (0,p]$ and any $s,t\in [0,T]$ we have
\begin{eqnarray*}
\E\left[|Y_t-Y_s|^r\right]&=&
\E\left[\Big|\int_s^t \zeta_\tau\,d\tau+\int_s^t \sigma(Y_\tau, \tau)dW_\tau\Big|^r\right]
\\
&\leq& C(r)\left\{\E\left[\Big|\int_s^t \zeta_t\,dt\Big|^r\right]+
\E\left[\Big|\int_s^t \sigma(Y_\tau, \tau)dW_\tau\Big|^r\right]
\right\}
\end{eqnarray*}
where $C(r)=2^{[r-1]_+}$ (notice that $r$ may be $<1$). Moreover, by the Burkholder-Davis-Gundy inequality and the bound on $\sigma$,
\begin{equation*}
\E\left[\Big|\int_s^t \sigma(Y_\tau, \tau)dW_\tau\Big|^r\right]
\le \E\left[\Big(\int_s^t\mbox{Tr}\big( \sigma(Y_\tau, \tau)\sigma^*(Y_\tau, \tau)\big)\,d\tau\Big)^{r\over 2}\right]
\le \delta^r\,|t-s|^{\frac r2}\,.
\end{equation*}
The conclusion follows combining the above estimates.
\hfill\QED

Let us finally recall some properties of Brownian bridges, which are one of the main ingredients of our method. 
\begin{Lemma}\label{BrownianBridge} Let $p\in(1,2)$ and let 
 $\sigma:\R^N\times [0,T]\to\R^{N\times D}$ be a Lipschitz continuous map such that $\|\sigma\|\le \delta$.
Then, for any $x, y \in \R^N$  there is 
a process $\zeta\in L^p_{\rm ad}(\Omega\times [0,T]; \R^N)$ such that the solution to 
$$
\left\{\begin{array}{l}
dY_t= \zeta_tdt+\sigma(Y_t, t)dW_t\\
Y_0=y
\end{array}\right.
$$
satisfies $Y_T=x$ $(\P\mbox{ a.s.})$ and
\begin{equation}\label{EstiZeta}
\E\left[\int_0^T |\zeta_t|^pdt\right]\leq C(p,\delta)\left(T^{1-p}|y-x|^p+ T^{1-p/2}\right).
\end{equation}
\end{Lemma}
Following \cite{FY}, $(Y_t)$ is called a Brownian bridge between $(y,0)$ and $(x,T)$. Estimate \eqref{EstiZeta}
can be found, e.g., in \cite{LS}. We give a proof of Lemma~\ref{BrownianBridge} for  completeness.

\medskip\noindent
{\it Proof:} Without loss of generality, we can assume that $x=0$. Having fixed $\alpha\in (1-1/p,2)$
(for instance $\alpha =3/4+1/(2p)$), let 
$Y_t$ be the solution to 
$$
\left\{\begin{array}{l}
dY_t= -\alpha\,\frac{Y_t}{T-t}\,dt+ \sigma(Y_t, t)dW_t\\
Y_0=y
\end{array}\right.
$$
We claim that 
\begin{equation}\label{FormuleY}
Y_t= T^{-\alpha}(T-t)^{\alpha}y+(T-t)^{\alpha}\int_0^t \sigma(Y_s,s)(T-s)^{-\alpha}dW_s
\end{equation}
and that (\ref{EstiZeta}) holds for
$
\zeta_t\doteq  -\,\alpha\,Y_t/(T-t)
$.
Indeed, let 
$$
Z_t=T^{-\alpha}(T-t)^{\alpha}y+(T-t)^{\alpha}\int_0^t \sigma(Y_s,s)(T-s)^{-\alpha}dW_s\;.
$$
Then $Z_0=y=Y_0$ and
\begin{eqnarray*}
dZ_t
&=& \Big(-\, \alpha \,T^{-\alpha}(T-t)^{\alpha-1}y -\, \alpha\,(T-t)^{\alpha-1}\int_0^t \sigma(Y_s,s)(T-s)^{-\alpha}dW_s\Big)dt 
 \\
& & \hspace{8cm}+\; (T-t)^{\alpha}\sigma(Y_t,t)(T-t)^{-\alpha}dW_t
 \\
& =&-\, \alpha\,(T-t)^{-1} Z_t\,dt+ \sigma(Y_t,t)\,dW_t 
\end{eqnarray*}
Hence, $Z_t=Y_t$ by uniqueness. Equality (\ref{FormuleY}) also implies that $Y_T=0$ ($\P$ a.s.). Let us now show that  (\ref{EstiZeta}) holds.  We have
$$
\zeta_t=  -\, \alpha\,\frac{Y_t}{T-t}= -\, \alpha \,T^{-\alpha}(T-t)^{\alpha-1}y-\, \alpha\,(T-t)^{\alpha-1} \int_0^t \sigma(Y_s,s)(T-s)^{-\alpha}dW_s
$$
Therefore,
\begin{eqnarray*}
\lefteqn{\E\left[\int_0^T |\zeta_t|^pdt\right]} \\
\; & \leq & 2^{p-1} \alpha ^pT^{-\alpha p}|y|^p \int_0^T (T-t)^{p(\alpha-1)}dt 
\\
& & \hspace{4.1cm}+2^{p-1} \alpha ^p\int_0^T (T-t)^{p(\alpha-1)}\E\left[\Big( \int_0^t \sigma(Y_s,s)(T-s)^{-\alpha}dW_s\Big)^p\right]dt\\
&\leq & C(p) T^{1-p}|y|^p 
\\
& & \hspace{.8cm}+ 2^{p-1}\,C(p)\, \alpha ^p\int_0^T (T-t)^{p(\alpha-1)} \E\left[\Big( \int_0^t {\rm Tr}(\sigma(Y_s,s)\sigma^*(Y_s,s))(T-s)^{-2 \alpha}ds\Big)^{p/2}\right]dt\\
&\leq & C(p) T^{1-p}|y|^p + C(p,\delta)T^{(1-2 \alpha)p/2}\int_0^T (T-t)^{p(\alpha-1)}dt\\
&\leq & C(p) T^{1-p}|y|^p + C(p,\delta)T^{1-p/2}\,,
\end{eqnarray*}
the second estimate above being justified by the Burkholder-Davis-Gundy inequality. 
\hfill\QED
%%%%%%%%%%%%%%%%%%%%%%%%%%%%%%%%%%%%%
%%%%%%%%%%%%%%%%%%%%%%%%%%%%%%%%%%%%%%
\section{Weak reverse H\"{o}lder inequalities}\label{se:holder}
Though sharing the same flavor of most results of common use,  the following reverse H\"older inequality lemma, obtained in \cite{ca08}, exhibits important differences in both assumptions and conclusion.  Since this is absolutely essential to our approach,  we will give a new proof of it which exploits a technique due to \cite{DaSb}. 
\begin{Lemma}\label{RevHolde} Let $p>1$ and let $\phi\in L^p(a,b)$ be a nonnegative function such that
\begin{equation}\label{eq:0hol_hyp}
 \frac{1}{t-a} \int_{a}^t \phi^p(s)\,ds \leq A\,\left( \frac{1}{t-a}\int_{a}^t \phi(s)\,ds\right)^p\qquad \forall t\in (a,b]
\end{equation}
for some constant $A>1$. Then, there are constants 
 $\theta=\theta(p,A)>p$ and $C=C(p,A)\geq 0$ such that
\begin{equation}
\label{eq:0hol_the}
\int_{a}^t\phi(s)\,ds \leq C\,(t-a)^{1-{1\over\theta}} (b-a)^{{1\over \theta}-{1\over p}}\,\|\phi\|_p\, \qquad \forall t\in [a,b]\;.
\end{equation}
\end{Lemma}

\begin{Remark}\rm\label{re:RevHolde} Observe that, by H\"older's inequality, 
\begin{equation*}
\int_{a}^t\phi(s)\,ds \leq (t-a)^{1-{1\over p}}\|\phi\|_p\, \qquad \forall t\in [a,b]\;.
\end{equation*}
So, the interest of the above lemma lies in the fact that \eqref{eq:0hol_the} provides the exponent
$1-1/\theta$ for $(t-a)$, which is higher than $1-1/p$.
\end{Remark}

\noindent {\it Proof: } Without loss of generality we can assume $a=0$ and $b=1$, the general form of the result being easy to recover by a rescaling argument. Let us further  assume that 
\begin{equation}\label{eq:extrahp}
\phi(t)\le \phi_0\qquad t\in[0,t_0]\quad\text{a.e.}
\end{equation}
for some constant $\phi_0\ge 0$ and some $t_0\in (0,1)$. 
Define
\begin{equation*}
f(s)={1 \over s}\int_0^s \phi(t)dt\qqf s\in (0,1]
\end{equation*}
and observe that, just like $\phi$, $f$ is bounded in a neighborhood of $0$.  
Now, let $\theta>p$ and recall Hardy's inequality (see, e.g.,\cite{hlp}) 
\begin{equation}
\label{Hardy}
\left( { \theta \over \theta-1}\right)^p\int_0^1 s^{{p\over\theta}-1}\ \phi^p(s)ds \geq  \int_0^1 s^{{p\over \theta}-1} \ f^p(s)ds\,.
\end{equation}
Moreover, observe that, in view of \eqref{eq:0hol_hyp}, 
\begin{multline}\label{eq:RHI}
 \int_0^1 s^{{p\over\theta}-1}f^p(s)ds \; \geq 
 \; {1\over A} \int_0^1 s^{{p\over\theta}-2}\left(\int_0^s\phi^p(t)dt\right)ds
 \\
 = {1\over A} \int_0^1 \phi^p(t) \left(\int_t^1 s^{{p\over\theta}-2}ds\right)dt 
 =  {\theta \over (p-\theta)A} \left( \int_0^1 \phi^p(t)dt - \int_0^1 t^{{p\over\theta}-1}\phi^p(t)dt \right).
\end{multline}
Then, combine \eqref{Hardy} and \eqref{eq:RHI}
to obtain
$$
\left[{\theta \over (p-\theta)A} +\left( { \theta \over \theta-1}\right)^p\right]\int_0^1 s^{{p\over\theta}-1}\phi^p(s)ds \; \geq \;
{\theta \over (p-\theta)A}   \int_0^1 \phi^p(s)ds
$$
or
$$
\left[{\theta \over (\theta -p)A} -\left( { \theta \over \theta-1}\right)^p\right]\int_0^1 s^{{p\over\theta}-1}\phi^p(s)ds \; \leq \;
{\theta \over (\theta-p)A}   \int_0^1 \phi^p(s)ds\,.
$$
Finally, choose $\theta=\theta(p,A)>p$ such that  
$${\theta \over (\theta -p)A} >\left( { \theta \over \theta-1}\right)^p$$
to deduce that, for some constant $C=C(p,A)$,
\begin{equation}\label{eq:step1}
\int_0^1s^{{p\over\theta}-1}\phi^p(s)ds \leq C\int_0^1 \phi^p(s)ds\,.
\end{equation}
At this point, the conclusion follows from H\"older's inequality and  
\eqref{eq:step1}: denoting by $q$ (resp. $\theta'$) the conjugate exponent of $p$ (resp. $\theta$), 
%i.e. $1/p+1/q=1$ (resp. $1/\theta+1/\theta'=1$), 
 we have
$$
\begin{array}{rl}
  \displaystyle{ \int_0^t \phi(s)ds} \; \leq &   \displaystyle{ \left( \int_0^t s^{p({1\over \theta}- {1\over p})}\phi^p(s)ds\right)^{1\over p} 
\left(\int_0^t s^{q({1\over p}-{1\over \theta})}ds\right)^{1\over q}  }\\
\leq & \displaystyle{   \left( \int_0^1 s^{{p\over \theta}- 1}\phi^p(s)ds\right)^{1\over p} \left({ \theta' \over q}\right)^{1 \over q} t^{1-{1\over \theta}} \; \leq 
  C^{{1\over p}}\left({ \theta' \over q }\right)^{1 \over q} t^{1-{1\over \theta}}\,\|\phi\|_p    }
\end{array}
$$
To complete the proof it remains to dispose of assumption \eqref{eq:extrahp}. 
For any $\tau\in (0,1]$, set
$$
\phi_\tau(s)=
\begin{cases}
\frac{1}{\tau}\int_0^\tau \phi(t)dt & {\rm if }\; s\in [0,\tau]
\vspace{2mm}
\\
\phi(s) & {\rm otherwise.}
\end{cases}
$$
Then, $\phi_\tau$ is bounded near 0 and $\phi_\tau\in L^p(0,1)$. We claim that  (\ref{eq:0hol_hyp}) is still true for $\phi_\tau$. Indeed, this is obvious if $t\in [0,\tau]$. On the other hand, for any
$t\in (\tau,1]$, 
\begin{multline*}
 \frac{1}{t}\int_0^t \phi^p_\tau(s)\,ds =\frac\tau t\left( \frac{1}{\tau}\int_{0}^\tau \phi(s)\,ds\right)^p+
 \frac{1}{t}\int_{\tau}^t \phi^p(s)\,ds   \leq 
 \frac{1}{t}\int_0^t \phi^p(s)\,ds
 \\
 \leq A\,\left( \frac{1}{t}\int_{0}^t \phi(s)\,ds\right)^p=A\,\left( \frac{1}{t}\int_{0}^t \phi_\tau(s)\,ds\right)^p\;.
\end{multline*}
Therefore, owing to the first part of the proof,
\begin{equation*}
\int_0^t \phi_\tau(s)ds \; \leq 
  C \,t^{1-{1\over \theta}}\, \|\phi_\tau\|_p   \qqf t\in [0,1]\,.
\end{equation*}
Letting $\tau\to 0^+$ gives \eqref{eq:0hol_the}.
\hfill\QED

We now give an adaptation of Lemma~\ref{RevHolde} that will be used in what follows.
\begin{Lemma}\label{RevHoldeBase} Let $p>1$ and let $\phi\in L^p(a,b)$ be a nonnegative function such that
\begin{equation}\label{eq:hol_hyp}
 \frac{1}{b-t} \int_{t}^b \phi^p(s)\,ds \leq A\,\left( \frac{1}{b-t}\int_{t}^b \phi(s)\,ds\right)^p+B\qquad \forall t\in [a,b)
\end{equation}
for some constants $A>1$ and $B\ge 0$. Then, there are constants\footnote{These are the same constants given by Lemma~\ref{RevHolde}.} 
 $\theta=\theta(p,A)>p$ and $C=C(p,A)\geq 0$ such that
\begin{equation}\label{eq:hol_the}
\int_{t}^b\phi(s)\,ds \leq C\,(b-t)^{1-{1\over\theta}}\,
\left\{(b-a)^{{1\over \theta}-{1\over p}}\|\phi\|_p+ B^{1/p}\,(b-a)^{1\over \theta}\right\}
\qquad \forall t\in [a,b]\,.
\end{equation}

\end{Lemma}
{\it Proof:} Let $\psi(s)\doteq \phi(s)+k$ where $k\doteq B^{1/p}/(A^{1/p}-1)$. In view of \eqref{eq:hol_hyp}, we have
\begin{eqnarray*}
\left( \frac{1}{b-t} \int_{t}^b \psi^p(s)\,ds\right)^{1/p}& \leq&\left( \frac{1}{b-t} \int_{t}^b \phi^p(s)\,ds\right)^{1/p} +k
\\
&\le & \frac{A^{1/p}}{b-t}\int_{t}^b \phi(s)\,ds+B^{1/p}+k
= \frac{A^{1/p}}{b-t}\int_{t}^b \psi(s)\,ds
\end{eqnarray*}
for every $t\in [a,b)$. Therefore, Lemma~\ref{RevHolde} applied to $\psi$ yields---after a change of variable---the existence of constants $\theta>p$ and $C\geq 0$, depending on $A$ and $p$ only, such that
\begin{equation*}
 \int_{t}^b \phi(s)\,ds = \int_{t}^b \psi(s)\,ds-k (b-t)\leq C\,(b-t)^{1-{1\over\theta}}
 \, (b-a)^{{1\over \theta}-{1\over p}}\,\|\psi\|_p\,.
\end{equation*}
The proof can now be completed noting that $\|\psi\|_p\le \|\phi\|_p+k(b-a)^{1/p}$.
\hfill \QED

We conclude this section with a generalization of Lemma \ref{RevHolde} to stochastic processes, which will be needed to study second order problems. 
\begin{Lemma}\label{lem:RevHol}
Let $(\Omega, {\cal A}, \P)$ be a probability space. Let $p\in (1,2)$ and let $\xi\in L^p(\Omega\times (a,b))$ be a nonnegative function such that
\begin{equation}\label{RevHolStoch}
\E\left[ \frac{1}{t-a} \int_{a}^t \xi_s^pds\right] \leq A\,\E\left[\Big(\frac{1}{t-a}\int_{a}^t\xi_sds\Big)^p\right]+\frac{B}{(t-a)^\frac{p}{2}}\qquad \forall t\in (a,b]
\end{equation}
for some positive constants $A$ and $B$. 
Then there are constants $\theta=\theta(p,A)\in(p,2)$ and $C=C(p,A)> 0$ such that
$$
\E \left[\Big(\int_{a}^t\xi_sds\Big)^p\right] \leq C (t-a)^{p-\frac{p}{\theta}}\left\{(b-a)^{\frac{p}{\theta}-1}\|\xi\|^p_p
+B(b-a)^{\frac{p}{\theta}-\frac{p}{2}}\right\}
\qquad \forall t\in (a,b]\;.
$$
\end{Lemma}
{\it Proof:} Under the extra assumption that $\xi_t$ is bounded ($\P$ a.s.) for a.e. $t$ near $a$, say for a.e. $t\in (a,t_0)$, let us define
\begin{equation*}
z_t={1 \over t-a}\int_a^t \xi_sds\qquad (\P\;\text{a.s.})\quad\forall t\in (a,b]\,.
\end{equation*}
Then, for any $\theta \in (p,2)$, Hardy's inequality \eqref{Hardy} yields
\begin{equation}
\label{Hardy_stoch}
\left( { \theta \over \theta-1}\right)^p\int_a^b (t-a)^{{p\over\theta}-1}\ \xi_t^p\,dt \geq  \int_a^b (t-a)^{{p\over \theta}-1} \ z_t^p\,dt \qquad (\P\;\text{a.s.})
\end{equation}
Owing to assumption \eqref{RevHolStoch}, we have
\begin{eqnarray}
\lefteqn{ \E\left[\int_a^b (t-a)^{{p\over \theta}-1} \ z_t^p\,dt\right]
  \; =\int_a^b (t-a)^{{p\over \theta}-1}\,
  \E\left[\Big(\frac{1}{t-a}\int_{a}^t\xi_sds\Big)^p\right]dt}
  \label{eq:RHI_stoch}
  \\\nonumber
&  \geq &
{1\over A} \int_a^b (t-a)^{{p\over \theta}-1}\left\{\E\left[ \frac{1}{t-a} \int_{a}^t \xi_s^pds\right]-\frac{B}{(t-a)^\frac{p}{2}}\right\}dt
 \\ \nonumber
& = &
{1\over A}\;\E\left[ \int_a^b (t-a)^{{p\over \theta}-2}\int_{a}^t \xi_s^p\,ds\, dt\right]-
 {B\over A} \,\int_a^b (t-a)^{{p\over \theta}-{p\over 2}-1} dt
 \\ \nonumber
& =&
{\theta\over (p-\theta)A}\;\E\left[ \int_a^b \left\{(b-a)^{{p\over \theta}-1} -(s-a)^{{p\over \theta}-1}\right\} \xi_s^p\,ds\right] -
{2\theta B\over p(2-\theta)A}  (b-a)^{{p\over \theta}-{p\over 2}} 
\end{eqnarray}
Now, taking the expectation of both sides of \eqref{Hardy_stoch} and chaining the resulting estimate with \eqref{eq:RHI_stoch} we obtain
\begin{multline*}
\left\{{\theta\over (p-\theta)A}+\left( { \theta \over \theta-1}\right)^p\right\} \E\left[\int_a^b (t-a)^{{p\over \theta}-1} \ \xi_t^p\,dt\right]
\\
\geq {\theta\over (p-\theta)A}\;(b-a)^{{p\over \theta}-1}\E\left[ \int_a^b  \xi_t^p\,dt\right] -
{2\theta B\over p(2-\theta)A}  (b-a)^{{p\over \theta}-{p\over 2}} 
\end{multline*}
or
\begin{multline*}
\left\{{\theta\over (\theta-p)A}-\left( { \theta \over \theta-1}\right)^p\right\} \E\left[\int_a^b (t-a)^{{p\over \theta}-1} \ \xi_t^p\,dt\right]
\\
\leq {\theta\over (\theta-p)A}\;(b-a)^{{p\over \theta}-1}\E\left[ \int_a^b  \xi_t^p\,dt\right] +
{2\theta B\over p(2-\theta)A}  (b-a)^{{p\over \theta}-{p\over 2}} 
\end{multline*}
Thus, choosing $\theta=\theta(p,A)\in(p,2)$ such that  
$${\theta \over (\theta -p)A} >\left( { \theta \over \theta-1}\right)^p$$
we deduce that, for some positive constant $C=C(p,A)$,
\begin{equation}\label{eq:step1_stoch}
 \E\left[\int_a^b (t-a)^{{p\over \theta}-1} \ \xi_t^p\,dt\right] 
 \leq C
 \left\{(b-a)^{{p\over \theta}-1}\E\left[ \int_a^b  \xi_t^p\,dt\right]
+B  (b-a)^{{p\over \theta}-{p\over 2}} \right\}\,.
\end{equation}
 By H\"older's inequality
and \eqref{eq:step1_stoch}, we have, for all $t\in (a,b]$,
\begin{multline*}
\E \left[\Big(\int_{a}^t\xi_sds\Big)^p\right]=
\E \left[\Big(\int_{a}^t(s-a)^{{1\over p}-{1\over \theta}}(s-a)^{{1\over \theta}-{1\over p}} \xi_sds\Big)^p\right]
\\ \le 
\E \left[\Big(\int_{a}^t(s-a)^{{q\over p}-{q\over \theta}}\Big)^{{p\over q}}
\int_{a}^t(s-a)^{{p\over \theta}-1} \xi_s^p\,ds\right]
\\
\le 
C(t-a)^{p-{p\over \theta}} \left\{(b-a)^{{p\over \theta}-1}\E\left[ \int_a^b  \xi_t^p\,dt\right]
+B  (b-a)^{{p\over \theta}-{p\over 2}} \right\}
\end{multline*}
which in turn implies the conclusion. Finally, our extra assumption that $\xi_t$ is bounded near $a$ can be removed arguing as in last part of the proof of Lemma~\ref{RevHolde}. \hfill\QED

%%%%%%%%%%%%%%%%%%%%%%%%%%%%%%%%%%%%%%%%%%%%%%
%%%%%%%%%%%%%%%%%%%%%%%%%%%%%%%%%%%%%%%%%%%%%%%%
\section{First order equations}\label{se:first}
In this section we shall be concerned with the first order Hamilton-Jacobi equation
\begin{equation}\label{HJ0}
u_t+H(x,t,Du)=0\qquad {\rm in }\quad \R^N\times (0,T)
\end{equation}
where 
\begin{equation*}
u_t={\partial u\over\partial t}\,,\quad Du=\Big({\partial u\over\partial x_1}\,,\,\dots\,,\,{\partial u\over\partial x_N}\Big)\,.
\end{equation*}
The data $H$ and $u$ will be assumed to satisfy the following hypotheses:
\begin{itemize}
\item $H:\R^N\times (0,T)\times \R^N\to\R$ is a continuous function such that 
\begin{equation}\label{GrowthCond}
\frac{1}{\delta}|z|^q-\eta_{_-} \leq H(x,t,z) \leq \delta |z|^q+\eta_{_+}\qquad \forall (x,t,z)\in \R^N\times (0,T)\times \R^N
\end{equation}
for some constants $q>1,\delta > 1$ and $\eta_\pm\ge 0$;
\item
$u:\R^N\times (0,T)\to\R$ is a bounded continuous viscosity solution of \eqref{HJ1}.
\end{itemize}
Moreover, we shall denote by  $p$ the conjugate exponent of $q$, i.e.,
\begin{equation*}
{1\over p}+{1\over q}=1\,,
\end{equation*}
and we shall fix a constant $M>0$ such that
\begin{equation}
\label{eq:M}
|u(x,t)|\le M\qquad\forall (x,t)\in\R^N\times (0,T)
\end{equation}
(e.g. $M=\sup_{(x,t)\in\R^N\times (0,T)}|u(x,t)|$). \\

In what follows, a {\em(universal) constant} is a positive number depending on the given data $q,\delta, N, \eta_\pm$ and $M$ only. Universal constants will be typically labeled with  $C$, but also with different letters (e.g., $\theta, A,\dots$). Dependence on extra quantities will be accounted for by using parentheses
(e.g., $C(r)$ denotes a constant depending also on $r$).

\begin{Theorem}\label{th:regu0}
Let  $u\in C(\R^N\times [0,T])$ be a viscosity solution of \eqref{HJ0}
satisfying \eqref{eq:M}. Then there is a universal constant  $\theta>p$ such that, for any $\tau\in (0,T]$, 
\begin{equation}\label{eq:main2}
|u(x,t)-u(y,s)|\leq C(\tau)\left[|x-y|^{\theta-p\over \theta-1}+|t-s|^{\theta-p\over\theta}\right]
\end{equation}
for all $(x,t), (y,s)\in \R^N\times[\tau,T]$ and some constant $C(\tau) >0$.
\end{Theorem}

\begin{Remark}\rm
The main point of the above theorem is that estimate \eqref{eq:main2} holds uniformly with respect to $H$ and $u$, as long as conditions \eqref{GrowthCond} and \eqref{eq:M} hold true. In particular, $C(\tau)$ is independent of the continuity modulus of $H$. 
\end{Remark}

%%%%%%%%%%%%%%%%%%%%%%%%%%%%%%%
\subsection{Proof of Theorem \ref{th:regu0}}

Let us set
\begin{equation*}
%\label{eq:H+}
H_{_+}(z)=\delta|z|^q+\eta_{_+}\qquad z\in\R^N
\end{equation*}
and
\begin{equation*}
%\label{eq:H-}
H_{_-}(z)=\frac{1}{\delta}|z|^q-\eta_{_-} \qquad z\in\R^N\,.
\end{equation*}
We begin the analysis with a kind of optimality  principle for super-solutions. 
\begin{Lemma}\label{le:1}
Let  $u\in C(\R^N\times[0,T])$ be a viscosity super-solution of 
\begin{equation}
\label{eq:1}
u_t+H_{_+}(Du)=0\qquad {\rm in }\quad \R^N\times(0,T)
\end{equation}
satisfying \eqref{eq:M}. Then, for any $(\bar x,\bar t)\in\R^N\times(0,T]$ there is an arc $\xi\in W^{1,p}([0,\bar t];\R^N)$, satisfying the end-point condition $\xi(\bar t)=\bar x$, such that
\begin{equation}\label{toto1}
u(\bar x,\bar t)\geq u(\xi(t),t) +C_{_+}\int_t^{\bar t}|\xi'(s)|^pds -\eta_{_+} (\bar t-t)\qquad \forall t\in [0,\bar t]
\end{equation}
where
\begin{equation}\label{eq:C+}
C_{_+}=\frac{\delta^{-p/q}}{pq^{p/q}}\,.
\end{equation}
\end{Lemma}
{\it Proof:} The requested arc will be provided by an approximation procedure.  
Let $(\bar x,\bar t)\in\R^N\times(0,T]$. For any positive integer $n$ let us set 
\begin{equation*}
\tau_n= {\bar t\over n}\,,\qquad t_k=k\tau_n\, \quad(k\in \{0,\dots, n\})\,.
\end{equation*}
%$$ (for $$) and $\tau=1/n$. 
We shall first construct a finite set of points $(x_k)_{k=0}^n$   
such that $x_n=\bar x$ and 
\begin{equation}
\label{constyk}
u(x_k,t_k)\geq u(x_{k-1},t_{k-1})+C_{_+}\tau_n\left|\frac{x_{k-1}-x_k}{\tau_n}\right|^p-\eta_{_+} \tau_n
\end{equation}
with $C_{_+}$ given by \eqref{eq:C+}. Having set $x_n=\bar x$,
we proceed as follows to construct $x_{k-1}$ from $x_k$ that we assume given for some $k\in \{1,\dots, n\}$. Let $v_k$ be the viscosity solution of
\begin{equation}\label{eq:pbmk}
\begin{cases}
v_t+H_{_+}(Dv)=0 & {\rm in }\quad \R^N\times(t_{k-1},T)
\\
v(x,t_{k-1})=u(x,t_{k-1}) & x\in\R^N\,.
\end{cases}
\end{equation}
As is well-known, $v_k$ is given by Hopf's formula
$$
v_k(x,t)=\min_{y\in \R^N}\left\{ (t-t_{k-1}) H_{_+}^*\left(\frac{y-x}{t-t_{k-1}}\right)+u(y,t_{k-1})\right\}\quad \forall (x,t)\in \R^N\times (t_{k-1}, T]\;,
$$
where $H_{_+}^*$ is the convex conjugate of $H_{_+}$, i.e.,
$$
H_{_+}^*(w)=\max_{z\in \R^N} \left\{ z\cdot w-H_{_+}(z)\right\}= C_{_+}|w|^p-\eta_{_+}
$$
with $C_{_+}$ given by \eqref{eq:C+}.
Since $u$ is a super-solution of \eqref{eq:pbmk}, the comparison principle yields $u(\cdot,t)\geq v_k(\cdot,t)$ for any $t\in[ t_{k-1},T]$. In particular, for $t=t_k$, we obtain that, for some
 point  $x_{k-1}\in\R^N$, 
$$
u(x_k,t_k)\geq v_k(x_k,t_{k})= u(x_{k-1},t_{k-1})+C_{_+}\tau_n\left|\frac{x_{k-1}-x_k}{\tau_n}\right|^p-\eta_{_+} \tau_n\,.
$$
The construction of $(x_k)_{k=0}^n$ can thus be completed by finite backward induction. 

Next, for any positive integer $n$,
let $\xi_n:[0,\bar t]\to \R^N$ be the piecewise linear interpolation of the above set
$(x_k)_{k=0}^n$ such that $\xi_n(t_k)=x_k$ for any $k=0,\,\dots\,,n$. We note that (\ref{constyk}) can be rewritten as
$$
u(\xi_n(t_k),t_k)\geq u(\xi_n(t_{k-1}),t_{k-1})+C_{_+}\int_{t_{k-1}}^{t_k} \left|\xi_n'(s)\right|^pds -\eta_{_+}(t_k-t_{k-1}) \quad (k=1,\,\dots\,,n)\,.
$$
Summing up the above inequalities gives
\begin{equation}\label{defxn}
u(\bar x,\bar t)\geq u(\xi_n(t_k),t_k)+C_{_+}\int_{t_k}^{\bar t} \left|\xi_n'(s)\right|^pds -\eta_{_+}(\bar t- t_k) \qquad \quad (k=0,\,\dots\,,n).
\end{equation}
Since $u$ is bounded,  $(\xi_n)_{n\in\N}$ is bounded in $W^{1,p}([0,\bar t];\R^N)$. Then, there is a subsequence of $(\xi_n)_{n\in\N}$ which weakly converges 
in $W^{1,p}([0,\bar t];\R^N)$ (hence converges uniformly) to some limit arc $\xi$ which satisfies $\xi(\bar t)=\bar x$. Passing to the limit in (\ref{defxn}) for such a subsequence gives (\ref{toto1}). 
\hfill \QED 

\begin{Remark}\rm\label{re:bound}
Observe that, owing to \eqref{eq:M}, for any arc $\xi$ which satisfies \eqref{toto1} we have
\begin{equation}\label{eq:Lp}
\int_{0}^{\bar t}|\xi'(s)|^pds\le 
{2M+\eta_{_+} \bar t\over C_{_+}}\,.
\end{equation}
\end{Remark}

We now turn to the analysis of sub-solutions of 
\begin{equation}
\label{eq:2}
u_t+H_{_-}(Du)=0\qquad {\rm in }\quad \R^N\times(0,T)
\end{equation}

\begin{Lemma}\label{lem:sub-sol}
Let $u\in C(\R^N\times[0,T])$ be a viscosity sub-solution of \eqref{eq:2}. 
Then, for any $s,t\in [0,T]$,  with $s> t$,
$$
u(x,s)\;\leq\; u(y,t)+C_{_-}(s-t)^{1-p}|y-x|^p+\eta_{_-}(s-t) \quad \forall x,y\in \R^N\;,
$$
where
\begin{equation}
\label{eq:C-}
C_{_-}=\frac{\delta^{p/q}}{pq^{p/q}}\,.
\end{equation}
\end{Lemma}
{\it Proof:} Let $\bar v$ be the viscosity solution of
\begin{equation}\label{eq:pbm-}
\begin{cases}
v_t+H_{_-}(Dv)=0 & {\rm in }\quad \R^N\times(t,T)
\\
v(x,t)=u(x,t) & x\in\R^N\,.
\end{cases}
\end{equation}
By Hopf's formula, 
\begin{equation}\label{eq:Hopf}
\bar v(x,s)=\min_{y\in \R^N}\left\{ (s-t) H_{_-}^*\left(\frac{y-x}{s-t}\right)+u(y,t)\right\}\quad \forall (x,s)\in \R^N\times (t, T)\,,
\end{equation}
where $H_{_-}^*$, the convex conjugate  of $H_{_-}$,  is given by
$$
H_{_-}^*(w)=\max_{z\in \R^N} \left\{ z\cdot w-H_{_-}(z)\right\}= C_{_-}|w|^p+\eta_{_-}\;.
$$
Since $u$ is a sub-solution of \eqref{eq:pbm-}, by comparison $u(\cdot,s)\leq \bar v(\cdot,s)$ for all $s\in [ t,T]$. 
So, \eqref{eq:Hopf} yields
$$
u(x,s)\leq \bar v(x,s) \leq C_{_-}(s-t)^{1-p}|y-x|^p+\eta_{_-}(s-t)+u(y,t)
$$
for all $(y,s)\in \R^N\times (t, T)$, as desired. \hfill\QED

Next, we derive a weak reverse H\"older inequality for the arcs that satisfy (\ref{toto1}).
\begin{Lemma}\label{Lemtoto2}
Let:
\begin{itemize}
\item $u\in C(\R^N\times[0,T])$ be a viscosity sub-solution of \eqref{eq:2} satisfying \eqref{eq:M};
\item
$\xi\in W^{1,p}([0,\bar t];\R^N)$  be an arc satisfying \eqref{toto1} with $\bar x=\xi(\bar t)$.
\end{itemize}
Then
\begin{equation}\label{toto2}
\frac{1}{\bar t-t} \int_t^{\bar t} |\xi'(s)|^pds \leq C_1\left(\frac{1}{\bar t-t}\int_t^{\bar t}|\xi'(s)|ds\right)^p+C_0
\qquad \forall t\in [0,\bar t)
\end{equation}
where 
\begin{equation}\label{def:C0C1}
C_0= \frac{\eta_{_+}+\eta_{_-}}{C_{_+}}   \quad \text{and}\quad C_1= \delta^{2p/q} >1\;.
\end{equation}
\end{Lemma}
{\it Proof:} Let $t\in [0, \bar t)$.  By Lemma \ref{lem:sub-sol}, 
$$
\begin{array}{rl}\displaystyle
u(\bar x,\bar t)\leq  &  u(\xi(t), t)+(\bar t-t)\,C_{_-}\left|\frac{\xi(t)-\bar x}{\bar t-t}\right|^p+\eta_{_-} (\bar t-t)\vspace{.5mm}\\
\leq  & \displaystyle 
u(\xi(t), t)+(\bar t-t)\,C_{_-}\left( \frac{1}{\bar t-t} \int_{t}^{\bar t} |\xi'(s)|ds \right)^p+\eta_{_-} (\bar t-t)\,.
\end{array}
$$
Combining (\ref{toto1}) with the above inequality we obtain
$$
C_{_+}\int_{t}^{\bar t}|\xi'(s)|^pds -\eta_{_+} (\bar t-t) \leq 
(\bar t-t) \,C_{_-}\left( \frac{1}{\bar t-t} \int_{t}^{\bar t} |\xi'(s)|ds \right)^p+\eta_{_-} (\bar t-t)\,,
$$
which in turns implies (\ref{toto2}). 
\hfill
\QED

In view of the above results, Lemma~\ref{RevHoldeBase} yields the following.
\begin{Corollary} \label{RevHolAppli}
Let  $u\in C(\R^N\times[0,T])$ be a viscosity solution of 
\begin{equation}
\label{eq:strip}
u_t+H(x,t,Du)=0\qquad {\rm in }\quad \R^N\times(0,T)
\end{equation}
satisfying \eqref{eq:M}, let $\tau\in (0,T)$ and let $(\bar x,\bar t)\in\R^N\times(\tau,T]$. Then  
there exist an arc $\xi\in W^{1,p}([0,\bar t];\R^N)$, satisfying the end-point condition $\xi(\bar t)=\bar x$ and inequality \eqref{toto1}, and 
a constant $\theta>p$, depending only on $q$ and $\delta$, such that 
\begin{equation}\label{co_toto2}
\int_t^{\bar t} |\xi'(s)|ds \leq C(\tau) \ (\bar t-t)^{1-{1\over\theta}}
\qquad \forall t\in [0,\bar t]
\end{equation}
for some constant  $C(\tau)\geq 0$.
 \end{Corollary}
{\it Proof:} First observe that, owing to 
\eqref{GrowthCond}, $u$ is a super-solution of \eqref{eq:1} and 
a sub-solution of \eqref{eq:2}. Then, Lemma~\ref{le:1} can be applied to construct 
an arc $\xi$ satisfying \eqref{toto1} together with  $\xi(\bar t)=\bar x$, while Lemma~\ref{Lemtoto2} ensures
that \eqref{toto2} holds true for $C_0,C_1$ given by \eqref{def:C0C1}. 
So, Lemma~\ref{RevHoldeBase} implies the existence of constant $\theta>p$ and $C_3>0$, depending only on
$C_1=\delta^{2p/q}$ and $p$,  such that 
\begin{equation*}
\int_t^{\bar t} |\xi'(s)|ds \leq C_3
\left\{\|\xi\|_p+k \bar t^{1\over p}\right\}
 {(\bar t-t)^{1-{1\over\theta}}\over \bar t^{{1\over p}-{1\over \theta}}}
\qquad \forall t\in [0,\bar t]\,,
\end{equation*}
where $k= C_0^{1/p}/(C_1^{1/ p}-1)$. 
Using the definition of $C_0, C_1$ and upper bound \eqref{eq:Lp} for $\|\xi\|_p$ yields \eqref{co_toto2}.
\hfill\QED

\noindent
{\it Proof of Theorem \ref{th:regu0} : } We will obtain two H\"older estimates in space and time, respectively, each of which will be uniform in the other variable. 

\medskip
\noindent
{\it Space regularity.}
Fix $\tau\in (0,T]$. Let $\bar t\in [\tau,T]$ and let $x,\bar x\in\R^N, x\neq\bar x$. 
   From Lemma~\ref{lem:sub-sol}, $u(x,\bar t)$ is bounded from above by
$$
u(x,\bar t)\;\leq\;u(y,t)+ C_{_-}(\bar t-t)^{1-p}|y-x|^p+\eta_{_-}(\bar t-t) \quad \forall y\in \R^N\;,
$$
for all $t\in [0,\bar t)$. Taking, in such an expression, $y= \xi(t)$, where $\xi$ is the arc provided by  the conclusion of Corollary~\ref{RevHolAppli}, yields, 
owing to \eqref{toto1}, 
$$
\begin{array}{rl}
u(x,\bar t)  \; \leq & \ds u(\xi(t),t) +  C_{_-} (\bar t-t)^{1-p} |\xi(t)-x|^p+\eta_{_-}(\bar t-t)
\vspace{.5mm}\\
\leq & u(\bar x,\bar t)+(\bar t-t)^{1-p} C_{_-}\left[\,|\xi(t)-\bar x|+|\bar x-x| \, \right]^p +\eta(\bar t-t)\end{array}
$$
for every $t\in[0,\bar t)$, where $\eta=\eta_{_-}+\eta_{_+}$. Therefore, since
$$
|\xi(t)-\bar x|\leq \int_t^{\bar t}|\xi'(s)|ds \leq 
C(\tau)\  (\bar t-t)^{1-\frac{1}{\theta}}
\qquad \forall t\in [0,\bar t]\;,
$$
on account of \eqref{co_toto2}, we obtain
\begin{equation}\label{eq:pre+xhol}
u(x,\bar t)-u(\bar x,\bar t)  \; \leq  \; (\bar t-t)^{1-p} C_{_-} \left\{ \,C(\tau)\ (\bar t-t)^{1-\frac{1}{\theta}}+|\bar x-x|  \, \right\}^p +\eta(\bar t-t)
\end{equation}
for every $t\in [0,\bar t]$.
Now, suppose $|\bar x-x|<\min\{1,\tau^{1-1/\theta}\}$. Then there is a number $t\in [0,\bar t]$ such that
$$\bar t-t=|\bar x -x|^{\theta\over\theta-1}\,.$$
Hence, \eqref{eq:pre+xhol} yields
$$
u(x,\bar t)-u(\bar x,\bar t)  \; \leq  \; C(\tau) |\bar x-x|^{\theta-p\over\theta-1}
%+2\delta |\bar x -x|^{\theta\over\theta-1}
%\qqf x\in\R^N
$$
for some new constant $C(\tau)>0$. On the other hand, the above inequality is trivial for 
$|\bar x-x|\ge \min\{1,\tau^{1-1/\theta}\}$ since $u$ is bounded. Moreover, the reasoning is symmetric with respect to $x$ and $\bar x$. So,
we have shown that
\begin{equation}\label{eq:xhol}
|u(x,\bar t)-u(\bar x,\bar t)|  \; \leq  \; C(\tau) |\bar x-x|^{\theta-p\over\theta-1}\qqf x,\bar x\in\R^N\,.
\end{equation}

\medskip
\noindent
{\it Time regularity.} Let $\bar x\in\R^N$ and let $\tau\leq t<\bar t\le T$.  Applying Lemma \ref{lem:sub-sol}
at $x=\bar x=y$, we obtain
\begin{equation}\label{eq:thol_0}
u(\bar x, \bar t)  - u(\bar x,t) \leq \eta_{_-} (\bar t-t)\,.
\end{equation}
To estimate the above left-hand side from below, let $\xi$ be as in the first part of the proof. Then,  
owing to \eqref{toto1},
\begin{equation}\label{eq:thol_1}
u(\xi(t),t)\leq u(\bar x,\bar t)+\eta_{_+}(\bar t-t)\,.
\end{equation}
On the other hand, in view of \eqref{eq:xhol} and \eqref{co_toto2},
\begin{equation}\label{eq:thol_2}
\begin{array}{rl}
u(\xi(t),t)\; \geq  & u(\bar x, t)- C(\tau) |\xi(t)-\bar x|^{\theta-p\over\theta-1}
\\
\ge & u(\bar x, t)- C(\tau)\left( \int_t^{\bar t} |\xi'(s)|ds\right)^{\theta-p\over\theta-1}
\;\ge \;u(\bar x, t)- C(\tau) (\bar t-t)^{\theta-p\over\theta}\,.
\end{array}
\end{equation} 
Combining \eqref{eq:thol_1} and \eqref{eq:thol_2} we conclude that
\begin{equation*}
u(\bar x, \bar t)  - u(\bar x,t) \geq -\eta_{_+} (\bar t-t)- C(\tau) (\bar t-t)^{\theta-p\over\theta}\,.
\end{equation*}
Since $(\theta-p)/\theta<1$, recalling \eqref{eq:thol_0} we finally get
\begin{equation}\label{eq:thol}
|u(\bar x,\bar t)-u(\bar x,t)|\leq C(\tau) |\bar t-t|^{(\theta-p)/\theta}
\qqf\bar x\in\R^N\,,\;\forall t,\bar t\in [\tau, T]\;.
\end{equation}
The conclusion follows from \eqref{eq:xhol} and \eqref{eq:thol}. \QED

%%%%%%%%%%%%%%%%%%%%%%%%%%%%%%%%%
\subsection{Local regularity for first order equations}
The H\"older regularity result of the previous section can be given a ``local version'', that is, a form that applies to solutions of the first order Hamilton-Jacobi equation
\begin{equation}\label{HJ1}
u_t+H(x,t,Du)=0\qquad {\rm in }\quad \O\,,
\end{equation}
where $\O$  is an open domain of $\R_x^N\times\R_t$. The Hamiltonian $H:{\cal O}\times \R^N\to \R$ is still assumed to satisfy the growth condition 
\begin{equation}\label{GrowthCondLoc}
\frac{1}{\delta}|z|^q-\eta_{_-} \;\leq\; H(x,t,z)\; \leq\; \delta |z|^q+\eta_{_+}\qquad \forall (x,t,z)\in \O\times \R^N
\end{equation}
for some constants $q>1,\delta \geq 1$ and $\eta_\pm\ge 0$.  Recall that $p$ is conjugate to $q$, and set
\begin{equation*}
\O_\rho=\{(x,t)\in\O~:~d_{\O^c}(x,t)>\rho\}\qqf \rho>0\,.
\end{equation*}

\begin{Theorem}\label{th:regu1} 
%Under assumption \eqref{GrowthCondLoc}, let  
Let $u$ be a continuous viscosity solution of \eqref{HJ1} 
satisfying $|u|\leq M$ in ${\cal O}$ for some $M\ge 0$.
Then there is a universal constant $\theta>p$ and, for any $\rho>0$, a constant $C=C(\rho)\ge 0$ such that
\begin{equation}\label{eq:main}
|u(x,t)-u(y,s)|\leq C\left[|x-y|^{\theta-p\over \theta-1}+|t-s|^{\theta-p\over\theta}\right]
\qqf (x,t), (y,s)\in \O_\rho \;.
\end{equation}
\end{Theorem}
{\it Proof:} Let $\rho>0$ be fixed and let  
$(\tilde x, \tilde t)\in {\cal O}_{4\rho}$. In order to simplify notation, we will assume that $\rho\leq 1/4$, $\tilde x=0$, and $\tilde t=2\rho$. Clearly, this implies no loss of generality. Set $T=4\rho$ and note that 
$$
B_{4\rho}\times (0,T)\subset \subset{\cal O}\;.
$$
Again without loss of generality, we can and will assume that the Hamiltonian has been extended to  $\R^N_x\times\R_t\times\R^N$, and that such an extension (still labeled by $H$) coincides with the original Hamiltonian on $B_{4\rho}\times (0,T)$ and satisfies \eqref{GrowthCond} on the whole space with the same constants that appear in \eqref{GrowthCondLoc}.

\medskip\noindent
{\sc  Step 1:} Let us show that there is a universal constant $\alpha>0$ such that, for any $(\bar x, \bar t)\in B_{2\rho}\times (0,T)$, there is an arc 
$\xi\in W^{1,p}([0,\bar t];\R^N)$, with $\xi(\bar t)=\bar x$, satisfying
\begin{equation}\label{Taillexi}
|\xi(t)-\bar x| \leq \rho
\end{equation}
and 
\begin{equation}\label{baruu}
u(\bar x,\bar t)\;\geq \;u(\xi(t),t) +C_{_+}\int_t^{\bar t}|\xi'(s)|^pds -\eta_{_+} (\bar t-t) 
\end{equation}
for all $t\in [(\bar t-\alpha  \rho^{p/(p-1)})_+,\bar t]$, where $C_{_+}$ is defined by  (\ref{eq:C+}). 

\medskip
\noindent
{\sc Proof:} Let $\phi(x)=\phi(|x|)$ be a smooth function defined on $\R^N$ such that $|\phi|\leq M$ in $\R^N$, $\phi\equiv M$ in $B_{3\rho}$, and $\phi\equiv-M$ in $B_{4\rho}^c$. Since $\phi$
is a trivial super-solution of (\ref{eq:1}), the function $\bar u:\R^N\times[0,T]\to\R$ defined by
$$
\bar u= 
\begin{cases}
 u\wedge \phi &\text{in}\quad B_{4\rho}\times [0,T]
 \\
 -M & \text{in}\quad B_{4\rho}^c\times [0,T]
\end{cases}
$$ 
is also a super-solution of (\ref{eq:1}) satisfying, thanks to (\ref{eq:M}), $\bar u\equiv u$ in $B_{3\rho}\times(0,T)$. So, applying Lemma \ref{le:1} to $\bar u$ we deduce that for any 
$(\bar x,\bar t) \in B_{2\rho}\times (0,T)$ there is an arc $\xi\in W^{1,p}([0,\bar t];\R^N)$, with $\xi(\bar t)=\bar x$,  such that
\begin{equation}\label{baru}
\bar u(\bar x,\bar t)\geq \bar u(\xi(t),t) +C_{_+}\int_t^{\bar t}|\xi'(s)|^pds -\eta_{_+} (\bar t-t)\qquad \forall t\in [0,\bar t]\;.
\end{equation}
Moreover, recalling Remark~\ref{re:bound},  
$$
\int_t^{\bar t}|\xi'(s)|^pds\leq 
{2M+\eta_{_+} \over C_{_+}}
\qqf t\in [0,\bar t]
$$
since $\rho\leq 1/4$. So, by H\"{o}lder's inequality,
\begin{equation}\label{eq:taillexi}
|\xi(t)-\bar x|\leq \left({2M+\eta_{_+} \over C_{_+}}\right)^{1\over p} (\bar t-t)^{1-{1\over p}}\qqf t\in [0,\bar t]\,.
\end{equation}
Therefore, taking 
$$
\alpha=\left({C_{_+}}\over 2M+\eta_{_+} \right)^{1\over p-1}\,,$$ 
inequalities (\ref{baru}) and \eqref{eq:taillexi}, combined with the fact that $\bar u=u$ in $B_{3\rho}\times (0,T)$, 
give \eqref{Taillexi} and \eqref{baruu}. 

Hereafter we will assume, without loss of generality,  that $\rho>0$ is such that 
$$\alpha\rho^{p/(p-1)} <\rho\,.$$

\medskip
\noindent
{\sc Step 2:} Let $(\bar x, \bar t)\in B_{2\rho}\times (\rho,T)$. We will show that  there exists universal constants $\theta>p$ and $C>0$ such that, if $\xi$ is an arc in $W^{1,p}([0,\bar t];\R^N)$ satisfying  \eqref{Taillexi}, \eqref{baruu} and $\xi(\bar t)=\bar x$ (as in Step 1), then
\begin{equation}\label{ineqxixi}
\int_t^{\bar t} |\xi'(s)|ds \;\leq\; C \, \rho^{-\, { \theta-p \over \theta(p-1) } }\, (\bar t-t)^{1-{1\over\theta}} 
\qquad \forall t\in [\bar t-\alpha \rho^{p/(p-1)},\bar t]\;.
\end{equation}
Moreover, for any $(x, t)\in B_{3\rho} \times[0, \bar t)$ and $y\in B_{3\rho}$,
\begin{equation}\label{turlututu}
u(x,\bar t) \;\leq  \;u(y,t)\,+\,C_{_-}\,(\bar t-t)^{1-p}\, |y-x|^p\,+\,C\,\rho^{-{p \over p-1}}\,(\bar t-t)\,.
\end{equation}
where $C_{_-}$ is defined by \eqref{eq:C-}.

\medskip
\noindent
{\sc Proof:} Let   $\phi$ be a function as in Step 1 such that $\|D\phi\|_\infty\leq C/\rho$ for some universal constant $C$. Then, 
$-\phi$ is a stationary sub-solution of 
\begin{equation}\label{eq:2tilde}
w_t+\widetilde H_{_-}(Dw)=0 \qquad {\rm in }\; \R^N\times (0,T)\;.
\end{equation}
where 
$$\begin{cases}
\widetilde H_{_-}(z)=\frac{1}{\delta}|z|^q-\widetilde \eta_{_-}& z\in\R^N
\vspace{1mm}
\\
\widetilde \eta_{_-}=\max\{\eta_{_-}, C^q/(\delta\rho^q)\}\,.
\end{cases}
$$
Let us set 
$$
\widetilde u=
\begin{cases}
u\wedge (-\phi) & {\rm in}\quad  B_{4\rho}\times (0,T)
\vspace{1mm}
\\
M & {\rm in} \quad  B_{4\rho}^c\times (0,T)\,.
\end{cases}
$$
Note that  $\widetilde u$ is a sub-solution of \eqref{eq:2tilde} such that $\widetilde u=u$ in  $B_{3\rho}\times (0,T)$, 
because $u$ is a sub-solution of (\ref{eq:2tilde}) and $|u|\le M$ in $B_{4\rho}\times (0,T)$. Let us now apply Lemma \ref{Lemtoto2} to $\widetilde u$,  $(\bar x,\bar t)$ and 
$\xi$: since $\widetilde u=u$ in $B_{3\rho}\times (0,T)$, 
$\xi$ satisfies 
$$
\widetilde u(\bar x,\bar t)\geq \widetilde u(\xi(t),t) +C_{_+}\int_t^{\bar t}|\xi'(s)|^pds -\eta_{_+} (\bar t-t)\qquad \forall t\in [\bar t-\alpha \rho^{p/(p-1)},\bar t]\,,
$$
we have
$$
\frac{1}{\bar t-t} \int_t^{\bar t} |\xi'(s)|^pds \leq C_1\left(\frac{1}{\bar t-t}\int_t^{\bar t}|\xi'(s)|ds\right)^p+C_0
\qquad \forall t\in [\bar t-\alpha  \rho^{p/(p-1)},\bar t)
$$
for some constants $C_0=(\eta_{_+}+\widetilde \eta_{_-})/C_{_+}=C_0'/ \rho^q$ and $C_1= \delta^{2p/q}>1$. 
Then, by Lemma \ref{RevHoldeBase}, we obtain the existence of  universal constants
$\theta>p$ and $C''>0$ such that 
$$
\int_t^{\bar t} |\xi'(s)|ds \leq C''\left(\|\xi\|_p+ k\alpha^{1\over p} \rho^{1 \over p-1}\right) \frac{(\bar t-t)^{1-{1\over\theta}}}{ \alpha^{ { \theta-p \over \theta p } } \rho^{ { \theta-p\over \theta(p-1) } } }
\qquad \forall t\in [\bar t-\alpha \rho^{p/(p-1)},\bar t]\;,
$$
where 
$$
k= \frac{C_0^{1\over p}}{C_1^{1\over p}-1}= \frac{C}{\rho^{q\over p}}=\frac{C}{\rho^{1\over p-1}} \quad
{\rm and }\quad \|\xi\|_p \le \left[ (2M+\eta_{_+})/C_{_+}\right]^{1\over p}\leq C
$$
for some universal constant $C$. Estimate (\ref{ineqxixi}) follows from the above inequality.  
Moreover, $\widetilde u$ being a sub-solution of (\ref{eq:2tilde}), Lemma \ref{lem:sub-sol} ensures that
$$
\widetilde u(x,\bar t) \; \leq \; \widetilde u(y,t)+(\bar t-t)^{1-p} C_{_-}|y-x|^p+\widetilde \eta_{_-}(\bar t-t)
$$
for any $(x, t)\in \R^N \times[0, \bar t)$, $y\in \R^N$. 
Since $\widetilde u=u$ in $B_{3\rho}\times (0,T)$, (\ref{turlututu}) follows for some constant $C$.

\medskip\noindent
{\bf Step 3:} We can now complete the proof of Theorem \ref{th:regu1}.

\smallskip
\noindent {\it Space regularity:}
Let $\bar t\in [\rho,T]$, let $x,\bar x\in B_{2\rho}$ be such that  $x\neq\bar x$, and let $\xi$ be the arc of Step 1. Taking $t\in[\bar t-\alpha \rho^{p/(p-1)},\bar t)$ and $y= \xi(t)$ in \eqref{turlututu}  yields 
$$
\begin{array}{rl}
u(x,\bar t)  \; \leq & \ds u(\xi(t),t)+ (\bar t-t)^{1-p} C_{_-} |\xi(t)-x|^p+C\rho^{-{p\over p-1}}(\bar t-t)
\vspace{1mm}\\
\leq & u(\bar x,\bar t)+ (\bar t-t)^{1-p} C_{_-}\left[\,|\xi(t)-\bar x|+|\bar x-x| \, \right]^p +C'\rho^{-{p\over p-1}} (\bar t-t) 
\end{array}
$$
for some universal constant $C'$. Hence, in view of \eqref{ineqxixi},
\begin{equation}\label{eq:pre+xhol2}
u(x,\bar t)-u(\bar x,\bar t)  \leq   (\bar t-t)^{1-p} C_{_-} \left\{ \,C \, \rho^{-\, { \theta-p \over \theta(p-1) } }\, (\bar t-t)^{1-{1\over\theta}} 
+|\bar x-x|  \, \right\}^p +C'\rho^{-{p\over p-1}} (\bar t-t)
\end{equation}
Now, suppose $|\bar x-x|<\alpha^{(\theta-1)/ \theta} \rho$. Then there is a number $t\in [\bar t-\alpha \rho^{p/ (p-1)}, \bar t]$ such that
$$\bar t-t=\rho^{ \theta-p \over (\theta-1)(p-1)}\; |\bar x -x|^{\theta\over\theta-1}\,.$$
So, owing to \eqref{eq:pre+xhol2},
$$
u(x,\bar t)-u(\bar x,\bar t)  \; \leq  \; C'' \left( \frac{|\bar x-x|}{\rho}\right)^{\theta-p\over\theta-1} +C'\left(\frac{|\bar x-x|}{\rho}\right)^{\theta\over \theta-1}
\leq C''' \left( \frac{|\bar x-x|}{\rho}\right)^{\theta-p\over\theta-1}
$$
for some new universal constants $C'', C'''>0$. Therefore,
$$
|u(x,\bar t)-u(\bar x,\bar t)|\leq C\left( \frac{|\bar x-x|}{\rho}\right)^{\theta-p\over\theta-1}
\quad \forall \bar x, x\in B_{2\rho}, \; \bar t\in [\rho, T]\;  {\rm  with }\; |\bar x-x|\le \alpha^{\theta-1 \over \theta}\rho \;.$$

\smallskip
\noindent {\it Time regularity:} Let $\bar x\in B_\rho$ and $\rho \leq t<\bar t\le T$.  Applying inequality \eqref{turlututu}
at $x=\bar x=y$, we obtain
$$
u(\bar x, \bar t)  - u(\bar x,t)\; \leq\; C\,\rho^{-{p\over p-1}} (\bar t-t)\,.
$$
To estimate the above left-hand side from below, let $\xi$ be given by Step 1. Then
$$
u(\xi(t),t)\;\leq\; u(\bar x,\bar t)+\eta_{_+}(\bar t-t)\,.
$$
Arguying as in the first step, we can choose a universal constant $\beta\in (0,\alpha)$
such that 
$$|\xi(t)-\bar x|\le \alpha^{\theta-1 \over \theta} \rho\qqf t\in [\bar t- \beta \rho^{p/ (p-1)},\bar t]\,.$$
Then, using the space regularity estimate we have just shown
and  \eqref{ineqxixi}, we obtain
\begin{multline*}
u(\xi(t),t)\;\geq\; u(\bar x, t)- C\left(\frac{|\xi(t)-\bar x|}{\rho}\right)^{\theta-p\over\theta-1}
\\
\qquad \ge\; u(\bar x, t)- C\left( \frac1\rho \int_t^{\bar t} |\xi'(s)|ds\right)^{\theta-p\over\theta-1}
\ge \;u(\bar x, t)- C\,\rho^{-{p(\theta-p) \over \theta(p-1)} }\,(\bar t-t)^{\theta-p\over\theta}\,.
\end{multline*}
Thus,
$$
|u(\bar x, \bar t)  - u(\bar x,t)| \;\leq C\;\,\rho^{-{p(\theta-p) \over \theta(p-1)} }\,(\bar t-t)^{\theta-p\over\theta}
$$
for all $\bar x\in B_\rho$ and all $\bar t, t \in [\rho, T]$ satisfying $|\bar t -t|\leq \beta \rho^{p/( p-1)}$. \hfill \QED

\begin{Remark}{\rm A simple analysis of the above proof allows to compute the dependence on $\rho$ of the constant in \eqref{eq:main} as follows
$$
|u(x,s)-u(y,t)|\leq C\left[ \rho^{-{\theta-p\over\theta-1}} |y-x|^{\theta-p\over\theta-1} +
\rho^{-{p(\theta-p) \over \theta(p-1) }} (t-s)^{\theta-p\over\theta}\right]
$$
for all $(x,s), (y,t) \in {\cal O}_\rho$ such that $|\bar x-x|\le k\rho $ and 
$|t-s|\leq k \rho^{p\over p-1}$, where  $C,k>0$ are universal constants.
}\end{Remark}

%%%%%%%%%%%%%%%%%%%%%%%%%%%%%%%%%
\section{Examples}\label{se:exa}
In this section we investigate two questions naturally arising from Theorems \ref{th:regu0} and \ref{th:regu1}. 
First, one may wonder whether the solutions  of \eqref{HJ1} satisfy stronger a priori estimates than \eqref{eq:main}, independent of the regularity of $H$. We address such a question with an example showing that uniform Lipschitz estimates cannot be expected even for a simple Hamilton-Jacobi equation in one space dimension.
Second, one may ask if  the local H\"{o}lder estimates for solutions in an open domain can be extended up to the boundary.
Surprisingly---and in stark contrast to the stationary setting (see \cite{CDLP})---this is not the case: we will exhibit a solution of a first order Hamilton-Jacobi
equation with constant coefficients which turns out to be discontinuous at the boundary of the domain.

\subsection{Counterexample to Lipschitz continuity}
\label{se:example_Lip}
The following example is inspired by \cite{AA}. In particular, Lemma \ref{le:ex1} and Proposition \ref{prop:ex1} could also be deduced from the results of the above paper. 

Let us fix $\gamma\in(2-\sqrt 2,1)$ and define
\begin{equation*}
\xi_0(t)=t^\gamma\qqf t\in[0,1]\,.
\end{equation*}

%\begin{figure}[htbp]
%\begin{center}
%{\includegraphics{goodarc}}
%\caption{$y=\xi_0(x)$}
%\label{fi:arc}
%\end{center}
%\end{figure}

\begin{Lemma}\label{le:ex1}
For every $t\in[0,1)$
\begin{equation}\label{eq:ex1}
\int_t^{t+h}|\xi'_0(s)|^2ds\;<\;{2\over h}\;|\xi_0(t+h)-\xi_0(t)|^2\qqf h\in(0,1-t]\,.
\end{equation}
\end{Lemma}
{\it Proof:} Let $t\in[0,1)$ and define
\begin{equation*}
X_t(h)=\int_t^{t+h}|\xi'_0(s)|^2ds\;-\;{2\over h}\;|\xi_0(t+h)-\xi_0(t)|^2\qqf h\in(0,1-t]\,.
\end{equation*}
Let us observe, first, that  $X_t(h)<0$ for $h>0$ small enough, since
$\lim_{h\downarrow 0}X_t(h)/h<0$. In order to obtain that $X_t(h)<0$ for every $h\in(0,1-t]$, let us show that $X_t(\cdot)$ is decreasing. Indeed, for any $h\in(0,1-t]$,
\begin{eqnarray}
X_t'(h)&=&|\xi'_0(t+h)|^2\;-\;{4\over h}\;\big(\xi_0(t+h)-\xi_0(t)\big)\xi'_0(t+h)
+\;{2\over h^2}\;|\xi_0(t+h)-\xi_0(t)|^2
\nonumber
\\
&=& \left(\frac{\xi_0(t+h)-\xi_0(t)}{h}\right)^2
\;\big[Y_t(h)^2-4Y_t(h)+2\big]
\label{eq:ex2}
\end{eqnarray}
where
\begin{equation*}
Y_t(h)\;\doteq\;\frac{h\,\xi'_0(t+h)}{\xi_0(t+h)-\xi_0(t)}.
\end{equation*}
Now, since $\xi_0$ is an incresing concave function, $Y_t(h)\le 1$. Moreover, 
\begin{equation*}
Y_t(h)\;\ge\; \frac{h\,\xi'_0(t+h)}{\xi_0(t+h)}\;=\;\frac{\gamma(t+h)^{\gamma-1}}{(t+h)^{\gamma-1}}\;=\;\gamma\;>\;2-\sqrt 2\,.
\end{equation*}
Since $y^2-4y+2<0$ for every $y\in (2-\sqrt 2,1]$, $X_t'(h)<0$ owing to \eqref{eq:ex2}.
\hfill\QED

Now, define
%\begin{equation*}
%L(x,t,v)=a(x,t)|v|^2\qqf (x,t,v)\in\R\times[0,1]\times\R
%\end{equation*}
%where
\begin{equation*}
a(x,t)=\begin{cases}
1 & \mbox{if}\quad x=\xi_0(t)
\\
2 & \mbox{if}\quad x\neq\xi_0(t)\,,
\end{cases}\qqf (x,t)\in\R\times[0,1]
\end{equation*}
and
\begin{equation*}
g(x)=\begin{cases}
0 & \mbox{if}\quad x=1
\\
G & \mbox{if}\quad x\neq 1
\end{cases}
\qqf x\in\R
\end{equation*}
where $G$ is a real number such that
\begin{equation}\label{eq:G}
G\;>\;{\gamma^2\over 2\gamma-1}\,.
\end{equation}
Let us consider the functional
\begin{equation*}
J[\xi] =\int_0^1a(\xi(t),t)|\xi'(t)|^2dt+g(\xi(1))\qqf \xi\in W^{1,2}([0,1])\,.
\end{equation*}

\begin{Proposition}\label{prop:ex1}
$\xi_0$ is the unique solution of the variational problem
\begin{equation}\label{eq:ex3}
\min\{ J[\xi]~:~\xi\in W^{1,2}([0,1])\,,\;\xi(0)=0\}\,.
\end{equation}
\end{Proposition}
{\it Proof:} To begin with, let us note that the minimum in \eqref{eq:ex3} does exists owing to well-known existence results for functionals with lower semicontinuous data (see, e.g., \cite[section 3.2]{BGH}). 
So, let $\xi_*$ be a solution of \eqref{eq:ex3} and observe that $\xi_*(1)=1$ since otherwise
\begin{equation*}
J[\xi_*]\ge g(\xi_*(1))=G\; >\;{\gamma^2\over 2\gamma-1}\;=\; J[\xi_0]\,.
\end{equation*}

Now, suppose that the open set $\{t\in(0,1)~:~\xi_*(t)\neq\xi_0(t)\}$ is nonempty and let $(t_1,t_2)$ be a connected component of such a set. Then, $\xi_*(t_i)=\xi_0(t_i)$ for $i=1,2$. Define
\begin{equation*}
\xi_1(t)\;=\;
\begin{cases}
\xi_*(t) &\mbox{if}\quad t\in[0,t_1]\cup [t_2,1]
\\
\xi_0(t) &\mbox{if} \quad t\in (t_1,t_2)\,.
\end{cases}
\end{equation*}
Then  $\xi_1\in W^{1,2}([0,1])$ satisfies $\xi_1(0)=0$ and $\xi_1(1)=1$. Moreover, in view of Lemma~\ref{le:ex1},
\begin{eqnarray*}
J[\xi_1]&=&\int_0^{t_1}a(\xi_*(t),t)|\xi_*'(t)|^2dt+\int_{t_1}^{t_2}\,|\xi_0'(t)|^2dt+\int_{t_2}^{1}a(\xi_*(t),t)|\xi_*'(t)|^2dt+g(1)
\\
&<&\int_0^{t_1}a(\xi_*(t),t)|\xi_*'(t)|^2dt+2\;{|\xi_0(t_2)-\xi_0(t_1)|^2\over t_2-t_1}+\int_{t_2}^{1}a(\xi_*(t),t)|\xi_*'(t)|^2dt+g(1)
\\
&=&\int_0^{t_1}a(\xi_*(t),t)|\xi_*'(t)|^2dt+2\;{|\xi_1(t_2)-\xi_1(t_1)|^2\over t_2-t_1}+\int_{t_2}^{1}a(\xi_*(t),t)|\xi_*'(t)|^2dt+g(1)
\\
&\le&J[\xi_*]
\end{eqnarray*}
in contrast with the optimality of $\xi_*$. Therefore, $\xi_*\equiv \xi_0$ and the proof is complete.
\hfill\QED

Let us now fix two sequences
\begin{equation*}
a_n:\R\times[0,1]\to\R\qmb{and}\qquad g_n:\R\to\R\qquad (n\ge 1)
\end{equation*}
of continuous functions such that
\begin{equation*}
\begin{cases}
{1\over 2} \le a_n(x,t)\le2&\forall n\ge 1
\\
a_n(x,t)\uparrow a(x,t) & n\to\infty
\end{cases}\qqf (x,t)\in\R\times[0,1]
\end{equation*}
and
\begin{equation*}
\begin{cases}
1\le g_n(x)\le g(x) &\forall n\ge 1
\\
g_n(x)\uparrow g(x) & n\to\infty
\end{cases}\qqf x\in\R\,.
\end{equation*}
For instance, one can take
\begin{equation}\label{eq:a_n}
a_n(x,t)=\min\Big\{2\;,\;n|x-\xi_0(t)|+\sum_{k=1}^n\frac1{2^k}\Big\}\qqf (x,t)\in\R\times[0,1]
\end{equation}
and
\begin{equation}\label{eq:g_n}
g_n(x)=\min\big\{G\;,\;n|x-1|\big\}\qqf x\in\R\,.
\end{equation}
Define, for all $n\in\N$,
\begin{equation*}
J_n[\xi] =\int_0^1a_n(\xi(t),t)|\xi'(t)|^2dt+g_n(\xi(1))\qqf \xi\in W^{1,2}([0,1])\,.
\end{equation*}
\begin{Proposition}\label{pr:exe1}
For every $n\in\N$ let $\xi_n$ be a solution of the variational problem
\begin{equation*}
\min\{ J_n[\xi]~:~\xi\in W^{1,2}([0,1])\,,\;\xi(0)=0\}\,.
\end{equation*}
Then $\xi_n\rightharpoonup \xi_0$ in $W^{1,2}([0,1])$ as $n\to\infty$.
\end{Proposition}
{\it Proof:} Since $(\xi_n)_n$ is bounded in $W^{1,2}([0,1])$, we can assume, without loss of generality, that $(\xi_n)_n$ weakly converges to some limit $\xi_*$ in $W^{1,2}([0,1])$. Consequently, $\xi_n\to\xi_*$ uniformly as $n\to\infty$.

Now, observe that,
\begin{equation}\label{eq:ex4}
\limsup_{n\to\infty}J_n[\xi_n]\le J[\xi_0]
\end{equation}
since $J_n[\xi_n]\le J_n[\xi_0]$ and, by monotone convergence,
$J_n[\xi_0]\to J[\xi_0]$ as $ n\to\infty$. Moreover, for any fixed $n\ge 1$ and all $m\ge n$,
$J_m[\xi_m]\ge J_n[\xi_m]$ in view of the monotonicity of $a_m$ and $g_m$. Therefore,
recalling \eqref{eq:ex4},
\begin{equation*}
J[\xi_0]\ge \liminf_{m\to\infty}J_m[\xi_m]\ge \liminf_{m\to\infty}J_n[\xi_m]\ge J_n[\xi_*]
\end{equation*}
owing to the lower semicontinuity of $J_n$. Since, by monotone convergence, $J_n[\xi_*]\to J[\xi_*]$ as $n\to\infty$, we conclude that $J[\xi_*]\le J[\xi_0]$. But we know that $\xi_0$ is the unique solution of \eqref{eq:ex3}. So, $\xi_*=\xi_0$ as requested.
\hfill\QED

Since $\xi_0$ is just H\"older continuous with exponent $\gamma$, and  $\xi_n\to \xi_0$ uniformly in $[0,1]$, the above result implies that $(\xi_n)_n$ cannot be equi-Lipschitz. 
\begin{Proposition}\label{pr:nolip}
Let $0<\tau<1$. Then the sequence of (value) functions
\begin{equation}\label{eq:un}
%\label{eq:exe5}
u_n(x,t)=\inf\left\{\int_t^1a_n(\xi(s),s)|\xi'(s)|^2ds+g_n(\xi(1))~:~\xi\in W^{1,2}([t,1])\,,\;\xi(t)=x\right\}
\end{equation}
is not equi-Lipschitz in $\R\times[0,\tau]$.
\end{Proposition}
{\it Proof:} Let $0<\tau<1$ and suppose $u_n$ is Lipschitz continuous in $\R\times[0,\tau]$, with the same constant
$K\ge 0$ for every $n\ge 1$. Let $\xi_n$ be as in Proposition~\ref{pr:exe1}. 
Then the optimality principle ensures that
$$u_n(0,0)=\int_0^ta_n(\xi_n(s),s)|\xi_n'(s)|^2ds+u(t,\xi_n(t))\qqf t\in [0,1]\,.$$ 
Therefore,
\begin{equation*}
\frac12\int_0^t|\xi_n'(s)|^2ds\le u_n(0,0)-u(t,\xi_n(t))\le K(t+|\xi_n(t)|)\qqf t\in [0,\tau]\,,
\end{equation*}
whence
\begin{equation*}
\left|\frac{\xi_n(t)}{t}\right|^2\le \frac 1t\int_0^t|\xi_n'(s)|^2ds\le 2K\left(1+\left|\frac{\xi_n(t)}{t}\right|\right)
\qqf t\in (0,\tau]\,.
\end{equation*}
The above inequality in turn implies that $|\xi_n(t)|\le c(K)t$ for every $t\in (0,\tau]$, uniformly for $n\ge1$, which is incompatible with the fact that $\xi_n\to \xi_0$ uniformly in $[0,1]$.
 \hfill\QED
 
Since $u_n$ above is the (unique) viscosity solution of the corresponding Hamilton-Jacobi equation, from Proposition~\ref{pr:nolip} we directly obtain the following corollary, which answers (negatively) the first question at the beginning of section~\ref{se:exa}. 
\begin{Corollary}\label{co:nolip}
For any integer $n\ge 1$ let $a_n$ and $g_n$ be given by \eqref{eq:a_n} and \eqref{eq:g_n}, respectively, and let $u_n$ be the viscosity solution of 
\begin{equation*}
\begin{cases}\ds
-u_t+{|u_x|^2\over 4a_n(x,t)}=0 &\mbox{in}\quad \R\times(0,1)
\vspace{.5mm}
\\
u(x,1)=g_n(x) & x\in\R\,.
\end{cases}
\end{equation*}
Then $(u_n)_n$ is not equi-Lipschitz in $\R\times[0,\tau]$, for any $0<\tau<1$.
\end{Corollary}
Observe that the above equation is of the form  \eqref{HJ1}, after the change of variable $t\mapsto 1-t$,  and satisfies  condition \eqref{GrowthCond} uniformly in $n$. 
\begin{Remark}\label{re:optimal_Holder}
\rm A careful examination of the proof of Proposition~\ref{pr:nolip} actually shows that
no uniform H\"older bound can be true for $(u_n)_n$ on $\R\times[0,\tau]$ with a H\"{o}lder exponent in the $x$ variable (resp. $t$ variable) greater than $1-1/\sqrt{2}$  (resp. 
$3-2\sqrt{2}$). Notice that such an optimal exponent is of the form $(\theta-2)/(\theta-1)$ 
(resp. $(\theta-2)/\theta$) for $\theta=1+\sqrt 2$ in agreement with \eqref{eq:main2}.

\end{Remark}

%%%%%%
\subsection{Counterexample to boundary continuity}
Our next example gives a negative reply to the second question raised at the beginning of section~\ref{se:exa}.
\begin{Example}\label{ex:bc}\rm
Let $\R_{_+}=(0,\infty)$ and consider the Hamilton-Jacobi equation
\begin{equation}\label{eq:bc}
u_t+\;\frac{1}4\;\Big|{\partial u\over\partial x}\Big|^2=0\qmb{in}\quad\O=\R_{_+}^2.
\end{equation}
Assumption \eqref{GrowthCond} is obviously satisfied with $q=2$. Now, define $u:\O\to\R$ by
\begin{equation*}
u(x,t)\doteq \min\Big\{1\;,\;\frac{x^2}{(t-1)_{_+}}\Big\}=
\begin{cases}
1 &\mbox{if}\;x^2\ge t-1
\\
\frac{x^2}{t-1} & \mbox{if}\;x^2< t-1
\end{cases}
\end{equation*}
Then, $u$ is a continuous function in $\O$ satisfying $0\le u\le 1$. Moreover, it is easily checked that $u$ is a solution of the above equation in $\O\setminus\Gamma$, where $\Gamma$ is the arc of parabola
\begin{equation*}
\Gamma=\{(x,t)\in \O~:~x^2= t-1\}\,.
\end{equation*}
So, since $u$ is locally semiconcave in $\O$ (see \cite{CannarsaSinestrari} for details), $u$ is a viscosity solution of \eqref{eq:bc}. On the other hand, $u$ is discontinuous at $(0,1)\in\partial \O$ because, for instance,
\begin{equation*}
\lim_{x\to 0^+}u(x,1)=1\qmb{while}\quad\lim_{x\to 0^+}u(x,1+2x^2)={1\over 2}\;.
\eqno{\Box}
\end{equation*}
\end{Example}
%%%%%%%%%%%%%%%%%%%%%%%%%%%%%%%%%

%%%%%%%%%%%%%%%%%%%%%%%%%%%%%%%%%
%%%%%%%%%%%%%%%%%%%%%%%%%%%%%%%%%
\section{Second order equations}\label{se:second}

%\subsection{Statement of the result}

In this section we are concerned with  second order Hamilton-Jacobi equations of the form 
\begin{equation}\label{HJ2}
u_t(x,t)-{\rm Tr}\left(a(x,t)D^2u(x,t)\right)+H(x,t,Du(x,t))=0\qquad {\rm in }\; \R^N\times (0,T)
\end{equation}
where 
$$
D^2u(x,t)=\left(\frac{\partial^2u}{\partial x_i\partial x_j}(x,t)\right)_{1\leq i,j\leq N}
$$
is the Hessian matrix. The data will be assumed to satisfy the following hypotheses:
\begin{itemize}
\item $H:\R^N\times (0,T)\times \R^N\to\R$ is a continuous function such that
\begin{equation}\label{GrowthCond2}
\frac{1}{\delta}|z|^q-\eta_{_-} \leq H(x,t,z) \leq \delta|z|^q+\eta_{_+}\qquad \forall (x,t,z)\in \R^N\times(0,T)\times \R^N\,,
\end{equation}
for some constants $q>2$, $\delta > 1$ 
and $\eta_\pm\geq 0$ ({\em super-quadratic growth});
\item there exists a locally Lipschitz continuous map $\sigma:\R^N\times (0,T)\to\R^{N\times D}$
 such that 
\begin{equation}\label{eq:cond_sigma}
a(x,t)=\sigma(x,t)\sigma^*(x,t)\quad\mbox{and}\quad
\|\sigma(x,t)\|\le \delta\qquad\forall (x,t)\in\R^N\times(0,T)\,;
\end{equation}
\item
$u:\R^N\times (0,T)\to\R$ is a continuous viscosity solution of \eqref{HJ2} such that 
$|u|\leq M$ in $\R^N\times[0,T]$.
\end{itemize}
As before,  {\it a universal constant} will be a positive number depending on the given data $q,\delta, M, \eta_{_-}, \eta_{_+}$ and $N$ only.  Recall that $p$ is the conjugate exponent of $q$. 

The main result of this section is the following H\"{o}lder estimate.
\begin{Theorem} \label{Regu2} Let $u\in C(\R^N\times[0,T])$ be a viscosity solution of \eqref{HJ2} such that 
$|u|\leq M$ in $\R^N\times[0,T]$. Then there is a universal constant $\theta>p$  
such that, for every $\tau>0$,
\begin{equation}\label{eq:main3}
|u(x_1,t_1)-u(x_2,t_2)|\leq C(\tau)\left[|x_1-x_2|^{(\theta-p)/(\theta-1)}+|t_1-t_2|^{(\theta-p)/\theta}\right]
\end{equation}
for any $(x_1,t_1), (x_2,t_2)\in \R^N\times [\tau, T]$ and for some constant $C(\tau) >0$.
\end{Theorem}
As for Theorem \ref{th:regu1}, the main point of the above result is  \eqref{eq:main3} holds true uniformly with respect to 
$H$ and $a$, as long as conditions \eqref{GrowthCond2} and the bound $|u|\leq M$ are  satisfied. In particular,  $\theta$ and  $C(\tau)$ are independent of the continuity moduli of $H$ and $a$. 

%%%%%%%%%%%%%%
\subsection{Some preliminary results}

For notational simplicity,  we prefer to replace the forward equation (\ref{HJ2})  by  the backward one
\begin{equation}\label{HJ3}
-u_t(x,t)-{\rm Tr}\left(a(x,t)D^2u(x,t)\right)+H(x,t,Du(x,t))=0\quad {\rm in }\quad \R^N\times (0,T)
\end{equation}
(which should be coupled with a {\it terminal condition}). Note that the change of variable $t\mapsto T-t$ turns a solution of (\ref{HJ2}) into a solution of (\ref{HJ3}), provided $a(x,t)$ and $H(x,t,z)$ are replaced by 
$a(x,T-t)$ and $H(x,T-t, z)$.

Throughout this section we shall need to keep track of the constants $\eta_{_+}$ and $\eta_{_-}$: indeed such a dependence is essential for the  proof of Theorem~\ref{ReguStochLoc}. 
For this purpose, we  will denote simply by $C$ (or $C_0$, $C_1$) constants which depend only on $\delta, M, p, T$ and $N$. Dependance with respect to
 $\tau$ and $\eta_{_\pm}$ will be made explicit by the use of parentheses. 

Let us begin with some estimates for super/sub--solutions of \eqref{HJ3}. 

\begin{Lemma}\label{lem:IneqStoch} Let $u\in C(\R^N\times [0,T])$ be a super-solution of 
\begin{equation}\label{eq:H+Stoch}
-u_t-{\rm Tr}\left(aD^2u\right)+\delta|Du|^q+\eta_{_+}=0\qquad {\rm in }\; \R^N\times (0,T)
\end{equation}
satisfying $|u|\leq M$ in $\R^N\times (0,T)$. Then, for any $(\bar x,\bar t)\in \R^N\times (0,T)$  there is a stochastic
basis $(\Omega,{\cal F}, \P)$, a filtration $({\cal F}_t)_{t\geq \bar t}$, a $D$-dimensional Brownian motion $(W_t)_{t\geq \bar t}$
adapted to the filtration $({\cal F}_t)$ and   a process  $\zeta\in L^p_{\rm ad}(\Omega\times [\bar t,T];\R^N)$, such that
the solution to
\begin{equation}\label{defX}
\left\{\begin{array}{l}
dX_t= \zeta_tdt+ \sqrt2\sigma(X_t,t)dW_t\\
X_{\bar t}=\bar x
\end{array}\right.
\end{equation}
satisfies
\begin{equation}\label{IneqStoch}
u(\bar x,\bar t)\geq \E\left[ u(X_t,t) +C_{_+}\int_{\bar t}^t|\zeta_s|^pds\right] -\eta_{_+} (t-\bar t)\qquad \forall t\in [\bar t, T]
\end{equation}
where $C_{_+}>0$ is the  universal constant given by \eqref{eq:C+}.
\end{Lemma}

\noindent {\it Proof: } Let $W$ be a $D$-dimensional Brownian motion on some probability space 
$(\Omega, {\cal A}, \P)$ and associated filtration $({\cal F}_t)$. Throughout the proof, for any
$y\in \R^N$, $t\in [0,T]$, $\zeta\in L^p_{\rm ad}(\Omega\times [t,T];\R^N)$ we denote by
$Y^{x, t, \zeta}$ the solution to
$$
\left\{\begin{array}{l}
dY_s= \zeta_sds+ \sqrt2\sigma(Y_s,s)dW_s\\
Y_{t}=x
\end{array}\right.
$$
Let $n$ be a large integer, 
$$
\tau=1/n\qquad {\rm and }\qquad t_k=\bar t+ k(T-\bar t)/n\qquad {\rm for}\; k\in \{0, \dots,n\}\;.
$$
Let us fix an initial condition $(\bar x,\bar t)\in\R^N\times [0,T]$. We are going to build 
a control $\zeta^n\in L^p_{\rm ad}(\Omega\times [\bar t,T];\R^N)$ and a process $Y^n= Y^{\bar x,\bar t, \zeta^n}$ such that
\begin{equation}\label{defYk}
u(Y^n_{t_k},t_k) \geq \E\left[ u(Y^n_{t_{k+1}},t_{k+1}) + C_{_+} \int_{t_k}^{t_{k+1}}\left|\zeta^n_s \right|^pds-(\eta_{_+}+\tau) \tau\; \big|\; {\cal F}_{t_k}\right]
\end{equation}
for any $k\leq n-1$. 

For any $k\geq 1$, let $v^k$ be the solution of (\ref{eq:H+Stoch}), defined on the time interval $[0,t_k]$, with terminal condition $u(\cdot, t_k)$. 
 From a classical representation formula (see, for instance, \cite{dale06}) we have
$$
v^k(x,t)=\inf_{\zeta\in L^p_{\rm ad}(\Omega\times [t, t_k];\R^N)} 
\E\left[ u(Y_{t_k}^{x,t, \zeta},t_k)+C_{_+}\int_{t}^{t_k}|\zeta_s|^pds-\eta_{_+} (t_k-t)\right]
$$
for any $(x,t)\in \R^N\times[0,t_k)$. 
Since $u(\cdot, t_k)$ is continuous, one can build, thanks to the measurable selection theorem (see \cite{aufr}), 
a Borel measurable map $x\to Z^{x,k}$ from $\R^N$ to $L^p_{\rm ad}(\Omega\times [t_{k-1}, t_k];\R^N)$
such that
$$
v^k(x,t_{k-1})\geq  \E\left[ u(Y_{t_k}^{x,t_{k-1}, Z^{x,k}},t_k)+C_{_+}\int_{t_{k-1}}^{t_k}|Z^{x,k}_s|^pds-(\eta_{_+} +\tau)\tau\right]\qquad \forall x\in \R^N\;.
$$
We now construct $Y^n$ and $\zeta^n$ by induction on the time intervals $[t_k, t_{k+1})$. On $[\bar t, t_1)$ we set
$\zeta^n_t=Z^{\bar x,\bar t}_t$ and  $Y^n= Y^{\bar x,\bar t,\zeta^n}$. Assume that $\zeta^n$ and $Y^n$
have been built on $[\bar t, t_{k-1})$. Then we set 
$$
\zeta^n= Z^{Y^n_{t_{k-1}},k} \qquad {\rm and }\qquad Y^n= Y^{\bar x,\bar t,\zeta^n}\qquad  {\rm on }\; [t_{k-1}, t_k)\;.
$$
Thus,  $(\zeta^n,Y^n)$ satisfies \eqref{defYk}. 

Next, from (\ref{defYk}) we get that 
\begin{equation}\label{Ynxin}
u(\bar x,\bar t) \geq \E\left[ u(Y^n_{t_{k}},t_{k}) + C_{_+} \int_{\bar t}^{t_{k}}\left|\zeta^n_s \right|^pds-(\eta_{_+}+\tau) (t_k-\bar t) \right]\qquad \forall k\leq n-1\;.
\end{equation}
In particular, since $u$ is bounded, we have the following bound for $\zeta^n$:
\begin{equation}\label{BorneZeta}
\E\left[\int_{\bar t}^T \left|\zeta^n_s \right|^pds\right]\leq C\qquad \forall n\geq 0
 \end{equation}
for some universal constant $C$. Let us set
$
\Lambda^n_t=\int_{\bar t}^t \zeta^n_sds
$ for all $t\in[\bar t, T]$.
Then \eqref{BorneZeta}, combined with Lemma \ref{stoch_bound}, implies that
$$
\E\left[\left|\Lambda^n_t-\Lambda^n_s \right|^p\right]\leq C|t-s|^{p-1}
 \qquad {\rm and }\qquad
\E\left[\left|Y^n_t-Y^n_s \right|^p\right]\leq C|t-s|^{p-1}
$$
for any $s,t\in [\bar t,T]$ and some universal constant $C$. 
Furthermore , since $Y^n_{\bar t}=\bar x$ and $\Lambda^n_{\bar t}=0$,  $(Y^n,\Lambda^n)$ satisfies Prokhorov's tightness condition. So, by 
Skorokhod's Theorem we can find a subsequence of 
$(\bar Y^n,\bar \Lambda^n)$ on some probability space $(\bar \Omega, \bar {\cal A}, {\bf \bar P})$
which has the same law as $(Y^n,\Lambda^n)$ and converges uniformly on $[\bar t,T]$ with probability 1 to some limit, say
$(\bar Y,\bar \Lambda)$. Since $\Lambda^n$ is absolutely continuous a.s. so is  
$\bar \Lambda^n$. Set $\bar \zeta^n_s= \frac{d}{ds} \bar \Lambda^n_s$. Then, 
 by \eqref{BorneZeta}, 
$
\E[\int_{\bar t}^T \left|\bar \zeta^n_s \right|^pds]\leq C
$ for all $n\geq 0$.
Therefore,  (up to a subsequence labeled in the same way) $(\bar \zeta^n)$ converges weakly in $L^p$ to some 
limit, $\bar \zeta$, which a.s. satisfies $\bar \Lambda_t=\int_{\bar t}^t\bar \zeta_sds$
for all $t\in [\bar t, T]$. Moreover, $\bar M^n_t:= \bar Y^n_t-\bar x-\bar \Lambda^n_t$
has the same law as $Y^n_t-\bar x-\Lambda^n_t=\int_{\bar t}^t\sigma(Y^n_s,s)dW_s$. Thus,  $\bar M^n$ is a continuous martingale satisfying
$$
\langle M^n_i,M^n_j\rangle_t= \int_{\bar t}^t (\sigma\sigma^*)_{ij}(\bar Y^n_s,s)ds\qquad \forall s\in [\bar t, T]\;.
$$
Hence, by \cite[Theorem 12]{mezh},  $\bar M:= \bar Y_t-\bar x-\bar \Lambda_t$ is also a continuous martingale with
$
\langle M_i,M_j\rangle_t= \int_{\bar t}^t (\sigma\sigma^*)_{ij}(\bar Y_s,s)ds
$ for all $s\in [\bar t, T]$.
Owing to the martingale representation theorem (see, e.g., \cite[Theorem 3.4.2]{KS}),
there is a $D-$dimensional Brownian motion $\widetilde W=\{\widetilde W_t, \widetilde {\cal F}_t, \ t\in [\bar t,T]\}$, 
defined on an extension $(\widetilde \Omega, \widetilde {\cal A}, {\bf \widetilde P})$ 
of $(\bar \Omega, \bar {\cal A}, {\bf \bar P})$, such that
$$
\bar Y_t-\bar x-\bar \Lambda_t= \bar M_t= \int_{\bar t}^t \sigma(\bar Y_s,s)d\widetilde W_s \qquad \forall t\in [\bar t, T]\;.
$$
Consequently, $
\bar Y_t=\bar x+\int_{\bar t}^t\bar \zeta_sds+\int_{\bar t}^t \sigma(\bar Y_s,s)d\widetilde W_s$ for all 
$t\in [\bar t, T]$. So, recalling \eqref{Ynxin}, a classical lower semicontinuity argument yields
$$
u(\bar x,\bar t) \geq \tilde \E\left[ u(\bar Y_{t},t) + C_{_+} \int_{\bar t}^{t}\left|\bar \zeta_s \right|^pds-\eta_{_+} (t-\bar t) \right]\qquad \forall t\in [\bar t,T]\;,
$$
which concludes the proof. \hfill \QED

\begin{Lemma}\label{lem:EstiSubSolStoch} Let $u\in C(\R^N\times [0,T])$ be a sub-solution of
\begin{equation}\label{eq:H-Stoch}
-u_t-{\rm Tr}\left(a D^2u\right)+\frac{1}{\delta}|Du|^q-\eta_{_-}=0\qquad {\rm in }\quad \R^N\times (0,T)
\end{equation}
satisfying $|u|\leq M$.
Then, 
for all $(\bar x,\bar t)\in \R^N\times (0,T)$  and all $(y,t)\in \R^N \times (\bar t, T)$,  
\begin{equation}\label{eq:EstiSubSolStoch}
u(\bar x,\bar t)  \leq  u(y,t)+C\left\{ |y-\bar x|^p(t-\bar t)^{1-p}+ (t-\bar t)^{1-p/2}\right\}+\eta_{_-}(t-\bar t)
\end{equation}
for some universal constant $C>0$.
\end{Lemma}

\begin{Remark}{\rm In particular, if $u$ is a sub-solution of the stationary equation
$$
-{\rm Tr}\left(a D^2u\right)+\frac{1}{\delta}|Du|^q-\delta=0\qquad {\rm in }\quad\R^N,
$$
then inequality (\ref{eq:EstiSubSolStoch}) implies that, for any $x,y\in \R^N$ and any  $\tau>0$, 
$$
u(x)  \leq  u(y)+C\left\{ |y-x|^p\tau^{1-p}+\tau^{1-p/2}\right\}+\eta_{_-}\tau\,,
$$
for some universal constant $C$. Thus, choosing $\tau=|x-y|^2$ yields
$
u(x)  \leq  u(y)+C\ |y-x|^{2-p}
$, that is, $u$ is H\"{o}lder continuous. This way we can partially recover one of the results in \cite{CDLP}. 
}\end{Remark}

\noindent {\it Proof:} Let us fix $t\in (\bar t,T)$. 
Let $v$ be the solution of equation (\ref{eq:H-Stoch}) with terminal condition $u(\cdot, t)$, and let  $(W_t)_{t\geq \bar t}$ 
be a $D$-dimensional Brownian motion on some stochastic basis $(\Omega, {\cal F}, \P)$, with
associated filtration $({\cal F}_t)_{t\geq \bar t}$. Then, by a classical representation formula (see, e.g., \cite{dale06}),
$$
v(\bar x,\bar t)= \inf_{\zeta\in L^p_{\rm ad}(\Omega\times [\bar t, t]; \R^N)} 
\E\left[ u(X_t^{\bar x,\bar t, \zeta},t)+C_{_-}\int_{\bar t}^t|\zeta_s|^pds+\eta_{_-} (t-\bar t)\right]
$$
where $C_{_-}$ is the constant given by (\ref{eq:C-}) and $X^{\bar x,\bar t, \zeta}$ is the solution of
\eqref{defX}. Owing to Lemma \ref{BrownianBridge}, we can choose $\zeta\in L^p_{\rm ad}(\Omega\times [\bar t, t]; \R^N)$ 
so that $X^{\bar x,\bar t, \zeta}$ is a Brownian bridge between $(\bar x,\bar t)$ and $(y, t)$ which satisfies
$$
\E\left[ \int_{\bar t}^t |\zeta_s|^p ds \right] \leq C\left\{ |y- \bar x|^p(t-\bar t)^{1-p}+ (t-\bar t)^{1-p/2}\right\}
$$
for some constant $C$ depending only on $p$ and $\delta$. Since $u$ is a sub-solution of \eqref{eq:H-Stoch}, the comparison principle yields
\begin{eqnarray*}
u(\bar x,\bar t)\leq v(\bar x,\bar t) & \leq & \E\left[ u(y,t)+C_{_-}\int_{\bar t}^t|\zeta_s(\omega)|^pds +\eta_{_-}(t-\bar t)\right]\\
&\leq & u(y,t)+C\left\{ |y-\bar x|^p(t-\bar t)^{1-p}+ (t-\bar t)^{1-p/2}\right\}+\eta_{_-}(t-\bar t)
\end{eqnarray*}
for some new constant $C$ (depending only on $p$ and $\delta$).
\hfill \QED

\begin{Lemma}\label{lem:EstiSubSolStoch2} Let $u\in C(\R^N\times [0,T])$ be a sub-solution of  \eqref{eq:H-Stoch} satisfying $|u|\leq M$. Fix
$(\bar x,\bar t)\in \R^N\times (0,T)$ and $\zeta\in L^p_{\rm ad}(\bar t, T;\R^N)$, and let $X$ be the solution  of \eqref{defX}. Then, for any $x\in \R^N$ and $t\in (\bar t, T)$,
\begin{equation}\label{eq:EstiSubSolStoch2}
\begin{array}{l}
u(x,\bar t)  - \E[u(X_t,t)]\\
\qquad \leq C\left\{(t-\bar t)^{1-p}\left(\E\left[\left(\int_{\bar t}^t |\zeta_s|ds\right)^p\right]+ 
 |\bar x-x|^p\right)+ (t-\bar t)^{1-p/2} \right\}
+\eta_{_-}(t-\bar t)
\end{array}
\end{equation}
for some constant $C>0$. 
\end{Lemma}
{\it Proof:}
Fix $t\in (\bar t,T)$ and apply Lemma \ref{lem:EstiSubSolStoch} to $y=X_t(\omega)$. Then, for almost all $\omega\in \Omega$,
$$
u(x,\bar t) \; \leq \; u(X_t(\omega),t)+C\left\{ |X_t(\omega)-x|^p(t-\bar t)^{1-p}+ (t-\bar t)^{1-p/2}\right\}+\eta_{_-}(t-\bar t)\;.
$$
Hence,
\begin{multline*}
\lefteqn{u(x,\bar t)}
\\
\leq \; \E\left[ u(X_t,t)\right]+
C\left\{ (\E\left[|X_t-\bar x|^p\right]+|\bar x-x|^p)(t-\bar t)^{1-p}+ (t-\bar t)^{1-p/2}\right\}+\eta_{_-}(t-\bar t).
\end{multline*}
Since, on account of Lemma \ref{stoch_bound},
$$
\E\left[|X_t-\bar x|^p\right]\leq C\left\{ \E\left[\left(\int_{\bar t}^t |\xi_s|ds\right)^p\right]+
\delta^p (t-\bar t)^{\frac{p}{2}} \right\},
$$
the conclusion follows. \hfill \QED

\begin{Lemma}\label{LemABStoch} Let $u\in C(\R^N\times [0,T])$ be a sub-solution of \eqref{eq:H-Stoch}
such that $|u|\leq M$ and let $\tau\in (0,T)$.
Then there is a universal constant $\theta\in(p,2)$ and a constant $C(\tau,\eta_{_\pm})>0$ such that, for every 
$(\bar x,\bar t)\in \R^N\times (0,T-\tau)$, every $\zeta\in L^p_{ad}(\bar t, T;\R^N)$, and every solution $X$ of \eqref{defX} satisfying \eqref{IneqStoch}, we have
$$
\E \left[\left(\int_{\bar t}^t|\zeta_s|ds\right)^p\right] 
\leq C(\tau,\eta_{_\pm})(t-\bar t)^{p-\frac{p}{\theta}}\qquad \forall t\in (\bar t, T]\;.
$$
\end{Lemma}
{\it Proof:} First, observe that, by Lemma \ref{lem:EstiSubSolStoch2} applied to $x=\bar x$, 
$$
u(\bar x,\bar t)  \leq  \E[u(X_t,t)]
+ C\left((t-\bar t)^{1-p}\E\left[\left(\int_{\bar t}^t |\zeta_s|ds\right)^p\right]+  (t-\bar t)^{1-p/2} \right)
+\eta_{_-}(t-\bar t)
$$
for all $t\in[\bar t,T]$. Moreover, in view of (\ref{IneqStoch}),
$$
\E\left[ u(X_t,t)\right] \leq u(\bar x,\bar t)-  C_{_+} \E\left[ \int_{\bar t}^t|\zeta_s|^pds\right] 
+\eta_{_+} (t-\bar t)\qquad \forall t\in [\bar t, T]\;.
$$
Hence,
$$
\E\left[ \int_{\bar t}^t|\zeta_s|^pds\right]  
\leq C_0\, (t-\bar t)^{1-p}\E\left[\left(\int_{\bar t}^t |\zeta_s|ds\right)^p\right]+ C_1(\eta_{_\pm})(t-\bar t)^{1-p/2}\qquad \forall t\in [\bar t, T]\;,
$$
for some constants $C_0$ and $C_1(\eta_{_\pm})$. Then, owing to Lemma \ref{lem:RevHol}, there
are constants $\theta\in(p,2)$ and $C(\eta_{\pm})>0$ such that 
$$
\E \left[\left(\int_{\bar t}^t|\zeta_s|ds\right)^p\right] \leq C(\eta_{\pm})\left(\|\zeta\|_p^p+1\right)
\frac{ (t-\bar t)^{p-\frac{p}{\theta}}}{(T-\bar t)^{1-\frac{p}{\theta}}}
\qquad \forall t\in (\bar t,T]\;.
$$
Since $u$ is bounded by $M$, assumption (\ref{IneqStoch}) implies that $\|\zeta\|_p\leq C(\eta_{_+})$. 
So, we finally get 
$$
\E \left[\left(\int_{\bar t}^t|\zeta_s|ds\right)^p\right] \leq C(\tau, \eta_{_\pm})(t-\bar t)^{p-\frac{p}{\theta}}\qquad \forall t\in (\bar t, T]\;,
$$
because $\bar t\leq T-\tau$.
\hfill \QED

%%%%%%%%%%%%%%%%%%%%%%%%%%%%%%%%%%%%%
\subsection{Proof of Theorem \ref{Regu2}} 

We are now ready to prove Theorem \ref{Regu2}. As above,  we will work with the backward equation
(\ref{HJ3}) instead of the forward one. Since $\eta_{_+}$ and  $\eta_{_-}$ are here fixed, we shall omit the dependence on such variables of all constants in  the proof. 

Let $u$ be a continuous, bounded solution of (\ref{HJ3}). Thanks to the growth assumption for $H$, $u$ is a super-solution of  (\ref{eq:H+Stoch}) and a sub-solution of (\ref{eq:H-Stoch}).
\\

\noindent {\it Space regularity:} Fix $(\bar x,\bar t)\in \R^N\times (0,T-\tau)$ and let $x\in \R^N$.  By Lemma \ref{lem:IneqStoch} 
there is a process $\zeta\in L^p_{\rm ad}(\Omega\times [\bar t,T];\R^N)$ and 
a solution $X$ to (\ref{defX}) such that (\ref{IneqStoch}) holds. So,
\begin{equation}\label{IneqStochSimple}
u(\bar x,\bar t)\geq \E\left[ u(X_t,t) \right] -\eta_{_+} (t-\bar t)\qquad \forall t\in [\bar t, T]\;.
\end{equation}
Also, Lemma  \ref{LemABStoch} ensures that 
\begin{equation}\label{EstiZeta4}
\E \left[\left(\int_{\bar t}^t|\zeta_s|ds\right)^p\right] \leq C(\tau)(t-\bar t)^{p-\frac{p}{\theta}}\qquad \forall t\in (\bar t, T]
\end{equation}
for some universal constant $\theta\in (p,2)$ and some constant $C(\tau)>0$. 
Furthermore,  applying Lemma~\ref{lem:EstiSubSolStoch2}, for any $t\in (\bar t, T]$ we have
$$
\begin{array}{l}
u(x,\bar t)  - \E[u(X_t,t)]\\
\qquad \leq C\left\{ (t-\bar t)^{1-p}\E\left[\left(\int_{\bar t}^t |\zeta_s|ds\right)^p\right]+ 
|\bar x-x|^p(t-\bar t)^{1-p}+ (t-\bar t)^{1-p/2}\right\}+\eta_{_-}(t-\bar t)\;.
\end{array}
$$
Plugging (\ref{IneqStochSimple}) and (\ref{EstiZeta4}) into the above inequality leads to
$$
u(x,\bar t) \leq u(\bar x,\bar t) +\eta(t-\bar t)+C(\tau)(t-\bar t)^{(\theta-p)/\theta}+C|\bar x-x|^p(t-\bar t)^{1-p}+C(t-\bar t)^{1-p/2}
$$
for any $t\in (\bar t, T)$, where $\eta=\eta_{_+}+\eta_{_-}$. Since $1>1-p/2>(\theta-p)/\theta$ (recall that $\theta<2$),
$$
u(x,\bar t) \leq u(\bar x,\bar t) +C(\tau)(t-\bar t)^{(\theta-p)/\theta}+C|\bar x-x|^p(t-\bar t)^{1-p}\;.
$$
Then, for $|x-\bar x|$ sufficiently small,  choose $t=\bar t+|x-\bar x|^{\theta/(\theta-1)}$ to obtain
$$
u(x,\bar t)\; \leq u(\bar x, \bar t)+C(\tau) |x-\bar x|^{(\theta-p)/(\theta-1)}.
$$

\noindent {\it Time regularity : } Let now $t\in (0,T-\tau)$. Then, in light of (\ref{IneqStochSimple}),
$$
 u(\bar x,\bar t) \geq \E\left[ u(X_t,t) \right]-\eta_{_+} (t-\bar t).
$$
Now, applying  the space regularity result we have just proved, we obtain
$$
\E\left[ u(X_t,t) \right]\geq  u(\bar x,t)- C(\tau) \E\left[ \left|X_t-\bar x\right|^{\theta-p\over\theta-1}\right].
$$
Moreover, by Lemma \ref{stoch_bound},
$$
\E\left[|X_t-\bar x|^{\theta-p\over\theta-1}\right] \leq  \E\left[\Big(\int_{\bar t}^t |\zeta_s|ds\Big)^{\theta-p\over\theta-1}\right]+ 
 C (t-\bar t)^{\theta-p\over 2(\theta-1)}\;.
$$ 
Also, by H\"{o}lder's inequality and (\ref{EstiZeta4}), 
$$
\E\left[\Big(\int_{\bar t}^t |\zeta_s|ds\Big)^{\theta-p\over\theta-1}\right]
\leq C\ \left\{\E\left[\Big(\int_{\bar t}^t |\zeta_s|ds\Big)^p\right]\right\}^{\frac{\theta-p}{p(\theta-1)}}
\leq 
C(\tau) (t-\bar t)^{\theta-p\over\theta}\;.
$$
Notice that $(\theta-p)/(2(\theta-1))> (\theta-p)/\theta$ since $\theta<2$. So,
$$
u(\bar x,\bar t)\geq u(\bar x,t) - C(\tau)  (t-\bar t)^{\theta-p\over\theta}\;.
$$
To derive the reverse inequality, one just needs to apply  Lemma~\ref{lem:EstiSubSolStoch} with $y=\bar x$ to get
$$
u(\bar x,\bar t)  \leq  u(\bar x,t) +C(t-\bar t)^{1-p/2}+\eta_{_-}(t-\bar t)\;.
$$
This leads to the desired result since $1-p/2>(\theta-p)/\theta$. 
\hfill \QED

%%%%%%%%%%%%%%%%%%%%%%%%%%
\subsection{Local regularity for second order equations}

We will now obtain a local version of Theorem \ref{Regu2}. Let ${\cal O}$ be a non-empty open subset of $\R^N_x\times \R_t$ and
consider the second order Hamilton-Jacobi equation 
\begin{equation}\label{HJ2bis}
u_t(x,t)-{\rm Tr}\left(a(x,t)D^2u(x,t)\right)+H(x,t,Du(x,t))=0\qquad {\rm in }\quad {\cal O}\;.
\end{equation}
As before $H:{\cal O}\times \R^N\to \R$ satisfies the super-quadratic growth condition
$$
\frac{1}{\delta}|z|^q-\eta_{_-} \leq H(x,t,z) \leq \delta|z|^q+\eta_{_+}\qquad \forall (x,t,z)\in {\cal O}\times \R^N,
$$
for some constants $q>2$, $\delta > 1$ 
and $\eta_\pm\geq 0$. Moreover, 
$a=\sigma\sigma^*$ where $\sigma:{\cal O} \to \R^{N\times D}$ (where $D\geq 1$) is assumed to be bounded by $\delta$
and locally Lipschitz continuous.

For any $\rho>0$, let us set
$
\O_\rho=\{(x,t)\in\O~:~d_{\O^c}(x,t)>\rho\}
$. 
\begin{Theorem} \label{ReguStochLoc} Let $u$ be a continuous viscosity solution of \eqref{HJ2bis} satisfying $|u|\leq M$
in ${\cal O}$ for some $M\geq 0$. Then there exists a universal constant  $\theta>p$ and, for any $\rho>0$, a constant $C(\rho)>0$
such that
$$
|u(x_1,t_1)-u(x_2,t_2)|\leq C(\rho)\left\{|x_1-x_2|^{\theta-p\over\theta-1}+|t_1-t_2|^{\theta-p\over\theta}\right\}
$$
for any $(x_1,t_1), (x_2,t_2)\in {\cal O}_\rho$. 
\end{Theorem}
{\it Proof:} Let us fix a point $(\tilde x, \tilde t)\in {\cal O}$ and let $\rho=\frac15d_{{\cal O}^c}(\tilde x, \tilde t)$. Without loss of generality we can assume that $\tilde x=0$ and $\tilde t =2T$, where $T=4\rho$. Notice that 
$$
\overline{B_{4\rho}}\times [0,T]\subset {\cal O}.
$$
Changing $t$ to $T-t$, we can also assume that $u$ is a solution of 
the backward equation
$$
-u_t(x,t)-{\rm Tr}\left(a(x,t)D^2u(x,t)\right)+H(x,t,Du(x,t))=0\qquad {\rm in }\quad B_{4\rho}\times (0,T).
$$
Without loss of generality we can extend $H$ and $\sigma$ outside of $\overline{B}_{4\rho}\times [0,T ]$ so that 
(\ref{GrowthCond2}) holds in $\R^N\times [0,T]$,
and $\sigma$ is bounded by $\delta$
and locally Lipschitz continuous in $\R^N\times [0,T]$.
\medskip

Let $\phi=\phi(|x|)$ be a smooth function such that 
$|\phi|\leq M$ in $\R^N$, $\phi=M$ in $B_{3\rho}$, and $\phi=-M$ in $B_{4\rho}^c$. 
Let $\bar \eta_{_+}(\rho)>\eta_{_+}$ and 
 $\tilde \eta_{_-}(\rho)>\eta_{_-}$ be such that 
$\phi$ is a (stationary) super-solution of  
\begin{equation}\label{eq:2bis}
-w_t(x,t)-{\rm Tr}\left(a(x,t)D^2w(x,t)\right)+\delta |Dw(x,t)|^q+\bar \eta_{_+}(\rho)=0\qquad {\rm in } \quad \R^N\times (0,T),
\end{equation}
while $-\phi$ 
 is a sub-solution of 
\begin{equation}\label{eq:2ter}
-w_t(x,t)-{\rm Tr}\left(a(x,t)D^2w(x,t)\right)+\frac{1}{\delta} |Dw(x,t)|^q-\tilde \eta_{_-}(\rho)=0\qquad {\rm in } \quad \R^N\times (0,T).
\end{equation}
Then, the map $\bar u$ defined by
$$
\bar u= \left\{\begin{array}{ll}
 u\wedge \phi & {\rm in }\; B_{4\rho}\times (0,T)\vspace{1mm}\\
-M & {\rm in }\;  B_{4\rho}^c\times (0,T)
\end{array}\right.
$$ is a super-solution of \eqref{eq:2bis} which satisfies $\bar u=u$ in $B_{3\rho}\times(0,T)$, whereas 
$$
\tilde u\doteq \left\{\begin{array}{ll}
 u\vee (-\phi) & {\rm in }\; B_{4\rho}\times (0,T)\vspace{1mm}\\
M & {\rm in }\;  B_{4\rho}^c\times (0,T)
\end{array}\right.
$$ is a sub-solution of \eqref{eq:2ter} such that $\tilde u=u$ in $B_{3\rho}\times(0,T)$.\medskip

Recall that, on account of Lemma \ref{lem:IneqStoch}, for every $(\bar x, \bar t)\in \R^N\times (0,T)$ there is a stochastic
basis $(\Omega,{\cal F}, \P)$, a filtration $({\cal F}_t)_{t\geq \bar t}$, a $D$-dimensional Brownian motion $(W_t)_{t\geq \bar t}$ 
 adapted to $({\cal F}_t)$,
and  a process  $\zeta\in L^p_{\rm ad}(\Omega\times [\bar t,T];\R^N)$ such that
the solution $X$ of (\ref{defX}) satisfies
\begin{equation}\label{IneqStochLoc}
\bar u(\bar x,\bar t)\geq \E\left[ \bar u(X_t,t) +C_{_+}\int_{\bar t}^t|\zeta_s|^pds\right] 
-\bar \eta_{_+}(\rho)(t-\bar t)\qquad \forall t\in [\bar t, T]
\end{equation}
with $C_{_+}$ given by \eqref{eq:C+}.

\medskip\noindent 
{\sc  Step 1:} Let $(\bar x, \bar t)\in B_{2\rho}\times (0,T)$  and let $X$ be as above. Then we claim that
\begin{equation}\label{eq:PX}
\P[|X_t-\bar x|\ge\rho] \leq C(\rho)(t-\bar t)^{p-1} \qquad \forall t\in [\bar t, T]
\end{equation}
(hereafter, $C(\rho)$ denotes a constant depending only on $N, q,\delta,M, \eta_\pm$ and $\rho$). In particular, 
\begin{equation}\label{barutildeu}
\left|\E[\bar u(t,X_t)]-\E[\tilde u(t,X_t)]\right|\leq  C(\rho) (t-\bar t)^{p-1}\qquad \forall t\in [\bar t, T]\;.
\end{equation}

\noindent {\sc Proof:}
Since $\bar u$ is bounded by $M$, (\ref{IneqStochLoc}) implies that 
\begin{equation}\label{boundzetap}
\E\left[ \int_{\bar t}^t|\zeta_s|^pds\right] \leq C(\rho) \qquad  \forall t\in [\bar t, T]\;.
\end{equation}
Now, by the Bienaym\'e-Chebyshev inequality and Lemma \ref{stoch_bound} we obtain
$$
\P[|X_t-\bar x|\geq \rho]\; \leq \; \frac{\E[|X_t-\bar x|^p]}{\rho^p}
\; \leq\; {C\over\rho^p} \left\{ \E\left[\left|\int_{\bar t}^t |\zeta_s|ds\right|^p\right]
+ (t-\bar t)^{p/2}\right\}.
$$
Also, by H\"{o}lder's inequality and (\ref{boundzetap}), 
$$
\E\left[\left|\int_{\bar t}^t |\zeta_s|ds\right|^p\right]\leq (t-\bar t)^{p-1} \E\left[\int_{\bar t}^t |\zeta_s|^pds\right]
\leq C(\rho)(t-\bar t)^{p-1}\;.
$$
Since $p/2>p-1$, \eqref{eq:PX} follows.
To show \eqref{barutildeu}, recall that $\bar u=\tilde u=u$ in $B_{3\rho}\times (0,T)$, with
$\bar u$ and $\tilde u$ bounded by $M$. Hence, owing to \eqref{eq:PX},
\begin{eqnarray*}
\left|\E[\bar u(t,X_t)]-\E[\tilde u(t,X_t)]\right| &\leq &
2M \P[X_t\notin B_{3\rho}] 
\\
&\leq& 2M \P[|X_t-\bar x|\ge \rho] \leq C(\rho) (t-\bar t)^{p-1}.
\end{eqnarray*}

\medskip
\noindent
{\sc Step 2:} Our next clam is that, for all $(\bar x,\bar t)\in B_{2\rho}\times (0,T-\rho)$ and 
$\zeta\in L^p_{\rm ad}(\Omega\times [\bar t,T];\R^N)$ such that
the solution $X$ to (\ref{defX}) satisfies \eqref{IneqStochLoc}, we have
\begin{equation}\label{ineqxixiStoch}
\E \left[\left(\int_{\bar t}^t|\zeta_s|ds\right)^p\right] \leq C(\rho) (t-\bar t)^{p-\frac{p}{\theta}}\qquad 
\forall t\in [\bar t, T]
\end{equation}
for some universal constant
$\theta\in(p,2)$. Moreover, for any $x\in B_{2\rho}$ and any $t\in (\bar t,T)$,
\begin{multline}\label{eq:EstiSubSolStoch3}
\tilde u(x,\bar t)  - \E[\tilde u(X_t,t)]\\
\leq C\left\{(t-\bar t)^{1-p}\left( \E\left[\Big(\int_{\bar t}^t |\zeta_s|ds\Big)^p\right]+ |\bar x-x|^p\right)+ 
(t-\bar t)^{1-p/2} \right\}+\tilde \eta_{_-}(t-\bar t)
\end{multline}

\noindent {\sc Proof:} 
Inequality (\ref{eq:EstiSubSolStoch3}) is a straightforward application of Lemma \ref{lem:EstiSubSolStoch2}.
Combining  (\ref{IneqStochLoc}) and (\ref{barutildeu}) gives
$$
\tilde u(\bar x,\bar t)\geq \E\left[ \tilde u(X_t,t) +C_{_+}\int_{\bar t}^t|\zeta_s|^pds\right] 
-\bar \eta_{_+}(\rho)(t-\bar t)-C(\rho)(t-\bar t)^{p-1}\qquad \forall t\in [\bar t, T]\;.
$$
Putting together the above estimate with (\ref{eq:EstiSubSolStoch3}), for $x=\bar x$, we obtain
$$
\E\left[\int_{\bar t}^t|\zeta_s|^pds\right]\leq 
C(t-\bar t)^{1-p} \E\left[\left(\int_{\bar t}^t |\zeta_s|ds\right)^p\right]+
C(\rho)(t-\bar t)^{1-p/2}\;.
$$
So, (\ref{ineqxixiStoch}) follows from Lemma \ref{lem:RevHol}, for some constants $\theta\in(p,2)$ and  $C(\rho)>0$. 

\medskip
\noindent
{\sc Step 3:} We can now proceed with the proof of Theorem~\ref{ReguStochLoc}.
Without loss of generality, in what follows we will assume that $\theta$ is so close to $p$ that
\begin{equation}\label{thetaclosep}
\frac{\theta -p}{\theta}\; \leq \;p-1\;.
\end{equation}

\noindent  {\it Space regularity:} 
Let $(\bar x,\bar t)\in B_{2\rho}\times (0,T-\rho)$ and let $x\in B_{2\rho}$.  Thanks to Step 1 
we can find $\zeta\in L^p_{\rm ad}(\Omega\times [\bar t,T];\R^N)$ and 
a solution $X$ to (\ref{defX}) satisfying (\ref{IneqStochLoc}). 
Since $\bar u=\tilde u=u$ in $B_{3\rho}\times (0,T)$, combining (\ref{IneqStochLoc}), (\ref{ineqxixiStoch}) and (\ref{eq:EstiSubSolStoch3}) with \eqref{barutildeu} yields, for any $t\in (\bar t, T)$,
\begin{multline*}
u(x,\bar t)  \\
 \leq u(\bar x,\bar t)+
C(\rho)\left\{ (t-\bar t)^{1-p}\left(\E\left[\left(\int_{\bar t}^t |\zeta_s|ds\right)^p\right]+ |\bar x-x|^p\right)+ 
(t-\bar t)^{1-p/2}+(t-\bar t)^{p-1}\right\}\\
\leq u(\bar x,\bar t)+
C(\rho)\left\{ (t-\bar t)^{1-p/\theta}+ |\bar x-x|^p(t-\bar t)^{1-p}\right\},
\end{multline*}
where we have also used (\ref{thetaclosep}) and the fact that $1-p/\theta<1-p/2$. 
For $t=\bar t+|x-\bar x|^{\theta/(\theta-1)}$ (for $|x-\bar x|$ sufficiently small) we then obtain
$$
u(x,\bar t)\; \leq u(\bar x, \bar t)+C(\rho) |x-\bar x|^{\theta-p\over\theta-1}\;.
$$
{\it Time regularity:} In view of (\ref{IneqStochLoc}) we have that
$$
u(\bar x,\bar t)\geq \E\left[ \bar u(X_t,t) \right]- \bar \eta_{_+} (t-\bar t)
$$
 for all $\bar t\in [0,T-\rho)$ and $t\in [\bar t,T]$.
Since $\bar u$ and $u$ are bounded functions that coincide on $B_{3\rho}$, recalling \eqref{eq:PX} we conclude that
\begin{eqnarray*}
\E\left[ \bar u(X_t,t) \right] &\geq &\E\left[ u(X_t,t){\bf 1}_{X_t\in B_{2\rho}} \right]
-M \P[|X_t-\bar x|\ge \rho] 
\\
&\geq& \E\left[ u(X_t,t){\bf 1}_{X_t\in B_{2\rho}} \right]
-C(\rho)(t-\bar t)^{p-1}\;.
\end{eqnarray*}
We now need to three further estimates. First, observe that, owing to our space regularity result above, 
\begin{eqnarray*}
\E\left[ u(X_t,t){\bf 1}_{X_t\in B_{2\rho}} \right] \; &\geq & 
\E\left[u(\bar x,t){\bf 1}_{X_t\in B_{2\rho}} \right]-C(\rho)\; \E\left[ |X_t-\bar x|^{\theta-p\over\theta-1}\right] \vspace{1mm}\\
 &\geq & u(\bar x,t) -C(\rho) \left\{ \E\left[ |X_t-\bar x|^{\theta-p\over\theta-1}\right]+ (t-\bar t)^{p-1}\right\}.
\end{eqnarray*}
Second, by Lemma \ref{stoch_bound}, 
$$
\E\left[|X_t-\bar x|^{\theta-p\over\theta-1}\right] \leq  C\;\E\left[\Big|\int_{\bar t}^t \zeta_s\,ds\Big|^{\theta-p\over\theta-1}\right]+ 
C (t-\bar t)^{\theta-p\over2(\theta-1)}.
$$ 
Third, by H\"older's inequality and  \eqref{ineqxixiStoch},
$$
\E\left[\Big(\int_{\bar t}^t |\zeta_s|ds\Big)^{\theta-p\over\theta-1}\right]\leq 
C(\rho) (t-\bar t)^{\theta-p\over\theta}.
$$
So, combining  all the above inequalities we obtain
$$
u(\bar x,\bar t)\geq u(\bar x,t) - C(\rho) \left\{(t-\bar t)^{p-1}+ (t-\bar t)^{\theta-p\over2(\theta-1)}+(t-\bar t)^{\theta-p\over\theta} \right\}- \bar \eta_{_+} (t-\bar t).
$$
Hence, recalling \eqref{thetaclosep} and  the fact that $(\theta-p)/(2(\theta-1))\geq (\theta-p)/\theta$, we get
$$
u(\bar x,\bar t)\geq u(\bar x,t) - C(\rho)  (t-\bar t)^{\theta-p\over\theta}.
$$

In order to show the reverse inequality, we just need to invoke
 Lemma~\ref{lem:EstiSubSolStoch}: indeed, taking $y=\bar x$ we get
$$
u(\bar x,\bar t)  \leq  u(\bar x,t) +C(t-\bar t)^{1-p/2}+\tilde \eta_{_-}(\rho)(t-\bar t)\;.
$$
The desired result follows since $1-p/2>(\theta-p)/\theta$.  
\hfill \QED
\\

\begin{center}
\bf Acknowledgements
\end{center}
The authors are greatly obliged to Catherine Rainer for continuous technical support. They are also grateful to Carlo Sbordone for bringing reference \cite{DaSb} to their attention.

%%%%%%%%%%%%%%%%%%%%%%%%%%%%%%%%%%%%%%%%%%%
%%%%%%%%%%%%%%%%%%%%%%%%%%%%%%%%%%%%%%%%%%
%\section{Appendix : Estimates for the Brownian bridges}

%{\bf POUR NOUS : }

\end{document}